\documentclass{article}
\usepackage{amssymb,latexsym} 
\newcommand{\1}{{{\mathchoice {\rm 1\mskip-4mu l} {\rm 1\mskip-4mu l} 
{\rm 1\mskip-4.5mu l} {\rm 1\mskip-5mu l}}}} 

\newcommand{\hol}{{holomorphic\,\,}}
\newcommand{\R}{{\mathbb R}} 
\newcommand{\PP}{{\mathbb P}}

\newcommand{\Q}{{\mathbb Q}}
\newcommand{\FH}{{\rm FH}} 
 
\newcommand{\SO}{{\rm SO}} 
\newcommand{\Flux}{{\rm Flux}} 
\newcommand{\Symp}{{\rm Symp}} 
 
\newcommand{\nnull}{{\rm null}} 

\newcommand{\coker}{{\rm coker}} 
 
\newcommand{\SU}{{\rm SU}}
\newcommand{\Gr}{{\rm Gr}}
\newcommand{\PSL}{{\rm PSL}}
\newcommand{\Jreg}{{\Jj_{{\rm reg}}}}
\newcommand{\GW}{{\rm GW}}
\newcommand{\QH}{{\rm QH}}
\newcommand{\univ}{{\rm univ}} 
\newcommand{\tu}{{\tilde{u}}}  
\newcommand{\tv}{{\tilde{v}}}
\newcommand{\sgrad}{{\rm s.grad}} 
 
\newcommand{\SL}{{\rm SL}} 
 
\newcommand{\tr}{{\rm tr}}

\newcommand{\Ad}{{\rm Ad}}
\newcommand{\LieHam}{{\rm LieHam}}
\newcommand{\rk}{{\rm rk}}

\newcommand{\bx}{{\bf x}}
\newcommand{\by}{{\bf y}}
\newcommand{\PU}{{\rm PU}}
\newcommand{\pbar}{{\ov {\p}}}
\newcommand{\ev}{{\rm ev}}

\newcommand{\TB}{{\widetilde B}}
\newcommand{\TJ}{{\widetilde J}}

\newcommand{\Z}{{\mathbb Z}} 
\newcommand{\C}{{\mathbb C}}

\newcommand{\ov}{\overline }

\newcommand{\G}{{\rm G}} 
\newcommand{\p}{{\partial}} 
\newcommand{\al}{{\alpha}} 
 
\newcommand{\be}{{\beta}} 
 
\newcommand{\Om}{{\Omega}} 
\newcommand{\om}{{\omega}}
\newcommand{\io}{{\iota}}
\newcommand{\eps}{{\varepsilon}} 
\newcommand{\de}{{\delta}} 
 
\newcommand{\ga}{{\gamma}} 
\newcommand{\Ga}{{\Gamma}} 
\newcommand{\ka}{{\kappa}} 
\newcommand{\la}{{\lambda}} 
\newcommand{\si}{{\sigma}}

\newcommand{\Dd}{{\mathcal D}} 
\newcommand{\Ll}{{\mathcal L}} 
\newcommand{\Jj}{{\mathcal J}}

\newcommand{\Nn}{{\mathcal N}} 
\newcommand{\Mm}{{\mathcal M}}

\newcommand{\Ss}{{\mathcal S}}

\newcommand{\La}{{\Lambda}} 
\newcommand{\Aut}{{\rm Aut}}

\newcommand{\Si}{{\Sigma}} 
\newcommand{\Ham}{{\rm Ham}} 
\newcommand{\area}{{\rm area\,}} 
\newcommand{\vol}{{\rm vol}}

\newcommand{\im}{{\rm Im }}

\newcommand{\IFF}{{\Longleftrightarrow}} 
\newcommand{\SSK}{{\smallskip}} 
\newcommand{\MS}{{\medskip}}

\newcommand{\NI}{{\noindent}} 
\newcommand{\proof}[1]{\noindent{\bf Proof#1:\ }} 
 
\newcommand{\QED}{\hfill$\Box$\medskip}

\newcommand{\Tilde}{\widetilde}
\newcommand{\THam}{{\Tilde{\Ham}}}
\newcommand{\TSymp}{{\widetilde{\Symp}}}

\newcommand{\trho}{{\Tilde{\rho}}}

\newcommand{\tphi}{{\Tilde{\phi}}}

\newcommand{\ta}{{\Tilde{a}}}

\newcommand{\Diff}{{\rm Diff}}

\newtheorem{theorem}{Theorem}[section] 
\newtheorem{cor}[theorem]{Corollary} 
 
\newtheorem{defn}[theorem]{Definition} 
\newtheorem{thm}[theorem]{Theorem}

\newtheorem{exercise}[theorem]{Exercise} 
 
\newtheorem{rmk}[theorem]{Remark}
\newtheorem{lemma}[theorem]{Lemma} 
\newtheorem{question}[theorem]{Question} 
 
\newtheorem{prop}[theorem]{Proposition}

\begin{document} 

\title{Lectures on groups of symplectomorphisms}
\author{Dusa McDuff\thanks{Partially
supported by NSF grant DMS 0072512}\\ State University of New York
at Stony Brook \\ (dusa@math.sunysb.edu)}

\date{July 18, 2003\\ revised version}

\maketitle
\MS\MS

\begin{abstract}
These notes combine material from short lecture courses 
given in  Paris, France, in 
July 2001 and in Srni, the Czech Republic, in January 2003.  
They discuss  groups 
of symplectomorphisms  of  closed symplectic manifolds $(M, \om)$
from various points of view.   
Lectures 1 and 2 provide an overview of our current knowledge 
of their algebraic, geometric and homotopy 
theoretic  properties.  Lecture 3 sketches the arguments
 used by Gromov, Abreu 
and Abreu--McDuff to figure out the rational homotopy type of these 
groups in the cases $M= \C P^2$ and $M=S^2\times S^2$. We outline the 
needed  $J$-holomorphic curve techniques.
Much of 
the recent progress in understanding the geometry and topology of 
these groups has come from studying the properties of fibrations with 
the manifold $M$ as fiber and structural group equal either to the symplectic 
group or to its Hamiltonian subgroup $\Ham(M)$. The case when the 
base is $S^2$ has proved particularly important.  Lecture 4 describes
the geometry of Hamiltonian fibrations over $S^2$, while Lecture 5 discusses 
their Gromov-Witten invariants via the Seidel representation.  It ends 
by sketching Entov's explanation of the ABW inequalities for 
eigenvalues of products of special unitary matrices.
Finally in Lecture 6 we apply the 
ideas developed in the previous two lectures to demonstrate the 
existence of (short) paths in $\Ham(M, \om)$ that minimize the Hofer norm over 
all paths with the given endpoints.  
\footnotetext{
\NI
keywords:  symplectomorphism group, 
Hamiltonian group, Hofer metric, symplectic fibrations, quantum 
homology, Seidel 
representation

}
\end{abstract}

\tableofcontents

\MS\MS

\NI
{\bf Acknowledgements}  The author warmly thanks the organisers of 
the schools in Jussieu and Srni for the wonderful 
working atmosphere  that they created.  She also thanks Kedra for 
some useful comments.
\MS

\section{Overview}

There are many different aspects to the study of groups of 
symplectomorphisms.  One can consider

\MS

\NI$\bullet$  algebraic properties; \SSK

\NI$\bullet$  geometric properties;\SSK

\NI$\bullet$  their homotopy type;\SSK

\NI$\bullet$  their stability under perturbations of the symplectic 
form.
\MS

One could also look at the dynamical properties of individual 
elements or of one parameter subgroups.
But we will not emphasize such 
questions here, instead concentrating on the above mentioned 
properties of the whole group.
In this article we will first survey some old and new results, 
mentioning some open problems, and then
will sketch some of the relevant proofs.  Background information 
and more references can be found in~\cite{MS0,MS1,MS2}.  Since the survey~\cite{Mcox}
discusses recent results on
the homotopy properties of the action of the Hamiltonian 
group on the underlying manifold $M$, this aspect of the theory will be only 
briefly mentioned here.

\subsection{Basic notions}

Throughout $(M,\om)$ will be a closed (ie compact and without boundary),
smooth (ie $C^\infty$),
symplectic manifold of dimension $2n$ 
unless it is explicitly mentioned otherwise.
The symplectomorphism group $\Symp(M,\om)$ 
consists of all 
diffeomorphisms $\phi: M\to M$ such that $\phi^*(\om) = \om$,
and is equipped with the $C^\infty$-topology.
Its identity component is denoted $\Symp_0(M,\om)$.  The latter 
contains an important  subgroup $\Ham(M,\om)$ whose elements are 
time $1$ maps of Hamiltonian flows $\phi_t^H$. 
(These are the flows $\phi_H^t, t\in [0,1],$ that at each time $t$ 
are tangent to
the symplectic gradient $X_t^H$ of the function $H_t:M\to \R$, i.e.
\begin{equation}\label{eq:H}
\frac \p{\p t} \phi_t^H (x) : = \dot{\phi}_t^H(\phi_t^H(x)) = 
X_t^H(\phi_t^H(x)), \qquad \om(X_t^H,\cdot) = -dH_t.)
\end{equation}
Thus we have a sequence of groups
and inclusions
$$
\Ham(M, \om)\hookrightarrow \Symp_0(M, \om)\hookrightarrow \Symp(M,\om)
\hookrightarrow \Diff(M).
$$
These groups are all infinite dimensional.  As explained in Milnor~\cite{Mil},
they can each be given the structure of a
Fr\'echet Lie group. As such they have well defined Lie 
algebras, with an exponential map.  For example, the Lie algebra of $\Ham(M, \om)$ consists of 
the space $C_0^{\infty}(M, \R)$  of smooth functions on $M$ with zero 
mean $\int_M H\om^n$, and the exponential map takes the (time 
independent) Hamiltonian  $H$ to the time $1$-map $\phi_1^H$ of the 
corresponding flow.  Observe that this map is never surjective.
\MS

The most important elementary theorems in symplectic geometry are:\SSK 

\NI
$\bullet$ {\bf 
Darboux's theorem}: {\it every symplectic form is locally 
diffeomorphic to the linear form 
$$
\om_0: = dx_1\wedge dx_2 + \dots + dx_{2n-1}\wedge dx_{2n}
$$
on Euclidean space}; and\SSK

\NI
$\bullet$ {\bf 
Moser's theorem}: {\it any path $\om_t, t\in [0,1],$ of {\it cohomologous} 
symplectic forms 
on a closed manifold $M$ is induced by 
an isotopy $\phi_t: M\to M$  of the underlying manifold, i.e. 
 $\phi_t^*(\om_t) = \om_0$, $\phi_0 = id$.  Forms $\om_0, \om_1$ that 
 are related in this way are called} {\bf isotopic}.
 \MS
 
 The fact that there are no local invariants of symplectic structures 
 is
 closely related to the fact that the symplectomorphism group is 
 infinite dimensional. Contrast this with Riemannian geometry in 
 which the curvature is a local invariant and isometry groups are 
 always finite dimensional.     Moser's theorem implies that 
the groups $\Symp(M,\om)$ and $\Ham(M, \om)$ depend 
only on the diffeomorphism class of the form $\om$.  In particular, 
they do not change their 
 topological or algebraic properties when $\om_t$ varies along a 
 path of cohomologous forms.  However, changes in 
 the cohomology class  $[\om]$ can cause 
 significant changes in the homotopy type 
 of these groups: see Proposition~\ref{prop:am2}. 
%  and the discussion of {\it robust} versus {\it fragile} in \S\ref{ss:act}.
In turn, this is closely related to the fact that the Moser theorem 
fails for families of noncohomologous forms.  As shown by 
McDuff~\cite{EX} (cf.~\cite[Chapter~9]{MS2}), there is a family of 
noncohomologous symplectic forms $\om_t, 0\le t\le 1,$ on $S^2\times 
S^2\times T^2$ such that $\om_0$ is cohomologous to $\om_1$ but not isotopic 
to it.
Here $\om_1$ is constructed to be diffeomorphic to $\om_0$, but the
construction can be modified to give  a family of 
 symplectic forms $\om_t, 0\le t\le 1,$ on  an $8$-dimensional 
 manifold  such that $\om_0$ is cohomologous to $\om_1$ but not 
 diffeomorphic to it.
 
\MS

Another basic point is that each of these groups have the homotopy type 
of a countable $CW$ complex (cf.~\cite[Chapter~9.5]{MS2}).

\subsection{Algebraic aspects}

Let $\TSymp_0(M)$ denote the universal cover of $\Symp_0(M)$.\footnote
{
We shall often 
drop $\om$ from the notation when it can be understood from the context.}
Its elements $\tphi$ are equivalence classes of paths $\{\phi_t\}_{t\in [0,1]}$
starting at the identity,
where $\{\phi_t\}\sim \{\phi_t'\}$ iff $\phi_1 = \phi_1'$ and the paths are 
homotopic with fixed endpoints.  We define
$$
\Flux (\tphi) = \int_0^1 [\om(\dot{\phi}_t,\cdot)] \in H^1(M, \R).
$$
One can check that $\Flux (\tphi)$
 is independent of the choice of 
representative for $\tphi$, and that the map 
$\tphi\mapsto\Flux(\tphi)$ defines a 
homomorphism $\TSymp_0(M) \to H^1(M, \R)$.
This is known as the {\bf 
Flux homomorphism}.   One way to see that $\Flux (\tphi)$ is well 
defined is to use the following alternate description.
Since $\Flux (\tphi)$ is a cohomology class,  it  is 
determined by its values on loops $s\mapsto\ga(s), s\in S^1,$ in $M$, 
and one can check that
\begin{equation}\label{eq:flux0} 
\Flux (\tphi)(\ga) = \int_{\tr_{\phi}(\ga)}\om,
\end{equation}
where $\tr_{\phi}(\ga)$ is the $2$-cycle 
$$
[0,1]\times S^1\to M: (t,s)\mapsto (\phi_t(\ga(s)).
$$
(For a proof of this and the other basic results in this 
section see~\cite[Chapter~10]{MS1}.)

One of the first results in the theory is that the rows and columns
in the following commutative diagram are short exact sequences of groups.
\begin{equation}\label{eq:flux}
\begin{array}{ccccc}
 \pi_1(\Ham(M))& \longrightarrow&
\pi_1(\Symp_0(M)) & \stackrel{\Flux}\longrightarrow &  \Ga_\om \\
 \downarrow& &\downarrow & & \downarrow  \\
 \THam(M)& \longrightarrow&
\TSymp_0(M) & \stackrel{\Flux}\longrightarrow&  H^1(M, \R)\\
\downarrow& &\downarrow & & \downarrow\\
\Ham(M)& \longrightarrow & \Symp_0(M) & 
\stackrel{\Flux}\longrightarrow& H^1(M, \R)/\Ga_\om.
\end{array}
\end{equation}
Here $\Ga_\om$ is the so-called {\bf flux group}.  It is the image of
$\pi_1(\Symp_0(M))$ under the flux homomorphism, and so far is not 
completely understood. In particular, it is not yet known whether 
$\Ga_{\om}$ is always discrete.  This question is discussed further 
in \S\ref{ss:stab}.  
\MS

One might wonder what other \lq\lq natural" homomorphisms there are
from $\Symp_0(M)$ to a arbitrary group $G$. 
If $M$ is closed, the somewhat 
surprising answer here is that {\it every} nontrivial homomorphism must factor 
through the flux homomorphism.  Equivalently, $\Ham(M)$ is {\it 
simple}, i.e. it has no proper normal 
subgroups.  The statement that $\Ham(M)$ has no proper closed normal 
subgroups is relatively easy and was proved by Calabi~\cite{CALAB}.
The statement that it has no proper normal subgroups at all is much more 
subtle and was proved by Banyaga~\cite{BANY}
 following a method introduced by Thurston 
to deal with the group of volume preserving diffeomorphisms.
The proof uses the relatively accessible fact\footnote
{
A general result due to Epstein~\cite{Eps} states that 
if a group $G$ of compactly supported homeomorphisms satisfies 
some natural axioms  then its 
commutator subgroup is simple.} 
that the commutator 
subgroup of $\Ham(M)$ is simple and the much deeper result that 
$\Ham(M)$ is a {\it perfect} group, i.e. is equal to its commutator subgroup.
More recently, Banyaga~\cite{BANY2} has shown 
that the manifold $(M, \om)$ may be 
recovered from the abstract discrete group $\Symp(M,\om)$. In other words,

\begin{prop}
If $\Phi: \Symp(M,\om) \to \Symp(M',\om')$
is a group homomorphism, then there is a diffeomorphism $f: M\to M'$
such that
$$
f^*(\om') = \pm \om,\qquad \Phi(\phi) = f\circ \phi\circ f^{-1} 
\;\;\mbox{ for  all }\phi\in \Symp(M,\om).
$$
\end{prop}

When $M$ is noncompact, the group $\Symp_0(M,\om)$
has many normal subgroups, for example the subgroup 
$\Symp_0(M,\om)\cap\Symp^c(M,\om)$ of 
all its compactly supported elements.  The identity component
$\Symp_0^c(M,\om)$ of the latter group supports a new homomorphism onto $\R$ 
called the Calabi homomorphism, and  Banyaga showed that the kernel of this 
homomorphism is again a simple group.  However, 
there is very little understanding of the normal 
subgroups of the full group
$\Symp(M)$.   In view of the above discussion of the
closed case the most obvious question is the following.

\begin{question} Is $\Symp(\R^{2n}, \om_0)$ a perfect
group?
\end{question}

It is known that the group of volume preserving 
diffeomorphisms of $\R^{k},k\ge 3,$ is perfect and is 
generated by elements whose 
support lies in a countable union of disjoint closed balls of radius 
$1$:
see McDuff~\cite{Mcvol} and Mascaro~\cite{Masc}.
(The same holds for the group of diffeomorphisms of $\R^n$ in any 
dimension.)  
However, the first of these statements is unknown in the symplectic case 
(even in the case of $\R^2$!)
while the second is false: see Barsamian~\cite{BARS}.

Some other questions of an algebraic nature are beginning to be tractable. 
Using ideas of Barge and Ghys, Entov~\cite{En2} has recently shown that 
in the closed case the universal cover of $\Ham(M, \om)$
has infinite commutator length, while  in~\cite{Pgd}
Polterovich develops methods to estimate the word length of 
iterates $f^{\circ n}$ in finitely generated subgroups of $\Ham(M,\om)$. 
This leds to interesting new restrictions on manifolds that support 
symplectic actions of nonamenable groups such as $\SL(2,\R)$. 

Entov and Polterovich~\cite{EnP} have also recently constructed a 
nontrivial continuous {\bf quasimorphism} $\mu$ on 
 $\Symp(S^2)$ and on the universal cover of $\Symp_0(M)$ for certain 
 other $M$. A quasimorphism on a 
group $\G$ is a map 
$\mu: \G\to \R$  that is a bounded distance away from being a 
homomorphism, i.e. there is a constant $c = c(\mu)> 0$ such that
$$
|\mu(gh) - \mu(g)-\mu(h)| < c, \qquad g,h\in \G.
$$
Their construction is particularly useful because $\mu$ restricts on 
to the Calabi homomorphism on the subgroups $G_U$ consisting of 
symplectomrophisms with support in a ball $U\subset M$. It
uses the structure of the quantum cohomology ring of 
$M$, and it is not yet clear whether it can be extended to all 
symplectic manifolds.  
%%!! added
Together with Biran, they recently showed~\cite{BEP} that 
such quasimorphisms give obstructions to displacing the Clifford 
torus in $\C P^n$.

\subsection{Geometric aspects}

The Lie algebra of the group $\Symp(M)$ is the space of all symplectic 
vector fields $X$, i.e. the vector fields on $M$ such that the $1$-form
$\io_X(\om)$ is  closed.  Similarly, the Lie algebra $\LieHam(M)$ is the 
space of all Hamiltonian vector fields $X$, i.e. those for which the $1$-form
$\io_X(\om)$ is  exact.   Since each exact $1$-form may be uniquely 
written as $dH$ where $H$ has zero mean $\int_M H \om^n$, $\LieHam(M)$ may 
also be identified with the space $C_0(M)$ of smooth functions on $M$ with 
zero mean.  With this interpretation, one easily sees that
$\LieHam(M)$ has a nondegenerate inner product
\begin{equation}\label{eq:inner}
\langle H,K\rangle = \int_M HK \om^n
\end{equation}
that is  bi-invariant under the adjoint action 
$$
\Ad_\phi(H) = H\circ \phi
$$
of $\phi\in \Symp(M)$.
Since a finite dimensional semisimple 
Lie group is compact if and only if it has such 
an inner product, this suggests that $\Ham(M)$ is an infinite dimensional 
analog of a compact group.   Of course this analogy is not perfect:
as we shall see in \S\ref{ss:rr} below, $\Ham(M)$ does not always 
have the homotopy type of a compact Lie group since it can have infinite 
cohomological dimension.  Nevertheless, it seems interesting to try to 
compare its properties with those of a compact Lie group.
  
One way in which  $\Ham(M)$ does resemble a compact Lie group is
that it supports a bi-invariant Finsler metric known 
as the {\bf Hofer metric}: see~\cite{Ho}. To define this, consider 
a path $\{\phi_t^H\}_{t\in [0,1]},$ in $\Ham(M)$
generated by 
the function $\{H_t\}_{t\in [0,1]}$.  Assuming that each $H_t$ has zero mean 
$\int_M 
H_t \omega^n$, we can define the negative and positive parts of its length 
by setting\footnote
{
To simplify the notation we will often write $H_t$ for the path 
$\{H_t\}_{t\in [0,1]}$ and $\phi_t$ instead of $\{\phi_t\}$.
Since it is seldom that we need to refer to $H_t$ for a fixed $t$ this should 
cause no confusion.}
$$
{\Ll}^-(H_t) = \int_0^1 -\min_{x\in M} H_t(x)\,dt,\quad
{\Ll}^+(H_t) = \int_0^1 \max_{x\in M} H_t(x)\,dt.
$$
Accordingly, we define seminorms $\rho^\pm$ and $\rho$
 by taking $\rho^\pm(\phi)$ 
to be the infimum of ${\Ll}^\pm(H_t)$
 over all Hamiltonians
$H_t$ with time $1$ map $\phi$, and  
 $\rho(\phi)$ to be the infimum  of 
$$
\Ll(H_t) = {\Ll}^+(H_t) + {\Ll}^-(H_t)
$$
over all such paths. 
It is easy to see that
$$
\rho^+(\phi) = \rho^-(\phi^{-1}),\quad
\rho^\pm(\phi\psi) \le \rho^\pm(\phi) + \rho^\pm(\psi),\quad
\rho^\pm(\psi^{-1}\phi\psi) = \rho^\pm(\phi).
$$
It follows that the metric $d_\rho(\phi,\psi): = \rho(\psi\phi^{-1})$
is bi-invariant and satisfies the triangle 
inequality.  Its nondegeneracy is equivalent to the statement
$$
\rho(\phi) = 0\;\;\; \IFF\;\;\; \phi = id.
$$
This deep result is the culmination of a series of papers by  
Hofer~\cite{Ho}, 
Polterovich~\cite{Pold} and Lalonde--McDuff~\cite{LMe}.  The key point is the following basic 
estimate which is known as the {\bf energy--capacity inequality}.

\begin{prop} Let $B = \phi(B^{2n}(r))$ be a symplectically
    embedded ball in $(M, \om)$ of  radius $r$.
    If $\phi(B)\cap B = \emptyset$, then  $\rho(\phi) \ge \pi r^{2}/2$.
\end{prop}

There has also been some success in describing the geodesics in 
$(\Ham(M), \rho)$.  This study was first initiated by 
Bialy--Polterovich, and a good theory has been developed for
paths that minimize length in their homotopy class.  
In \S\ref{sec:lect6} we shall sketch the proof that absolutely length 
minimizing paths exist.  Here is a simple form of the result.  (It has 
been recently generalised by Oh~\cite{Oh4} to  {\bf quasiautonomous} 
Hamiltonians.)

\begin{prop} The 
natural $1$-parameter subgroups $\{\phi_t^H\}_{t\in \R}$
generated by time independent $H$ minimize length between the identity 
and $\phi_t^H$ for all sufficiently small $|t|$.  Thus for each $H$ there is $T= 
T(H)>0$
such that
$$
\rho(\phi_t^H) = \Ll(\{H\}_{t'\in [0,|t|]}) = |t|(\max H - \min H)
$$
whenever $|t|\le T$.
\end{prop}

There are still many interesting open questions about 
Hofer geometry, some of which are mentioned below. 
Interested readers should consult  Polterovich's book~\cite{Pbk} for 
references and further discussion.  There are also beginning to be 
very interesting dynamical applications of Hofer geometry: see, for 
example, Biran--Polterovich--Salamon~\cite{BPS}.

\begin{question} Does $\Ham(M)$ always have 
infinite diameter with respect to 
the Hofer norm $\rho$?
\end{question}

This  basic question  has not yet been answered 
because of the difficulty of finding lower bounds for $\rho$.  The 
most substantial result here
is due to Polterovich, who gave an affirmative answer in the case
 $M = S^2$: see~\cite{Pdm,Pbk}.  
\MS

\begin{question}\label{q:slim}
Find ways to estimate the maximum value of 
$T$ such that the path $\{\phi_t^H\}_{t\in [0,T]}$
minimizes the length between the identity and $\phi_T^H$ among all 
homotopic paths  with  fixed endpoints.
\end{question}

If $H$ is time independent it has been shown by 
McDuff--Slimowitz~\cite{MSlim} and Entov~\cite{En} 
that one can let $T$ be the smallest positive number such that either 
the flow $\phi_t^H$ of $H$ or one of the linearized flows at its critical points
 has a nontrivial periodic orbit of period $T$.
 A Hamiltonian that satisfies this condition with $T=1$ is said to be 
 {\bf slow}.  
 If a compact Lie group $\G$ acts  effectively
 on $(M, \om)$ by Hamiltonian 
symplectomorphisms then its image in $\Ham(M)$ consists of elements 
that are the time-$1$ maps of the flows of autonomous Hamiltonians.
Hence its image is totally geodesic 
with respect to the Hofer norm.  Moreover the restriction of the Hofer 
norm to $\G$ often has interesting geometric properties.  For an 
example, see the discussion of the ABW inequalities in Section~\ref{ss:ABW}.
 
Oh has made considerable progress with the general question (for 
nonautonomous flows)
 in his recent paper\cite{Oh4}.  He has also defined a refined version of the 
 Hofer norm using spectral invariants, that coincides with it in a 
 neighbourhood of the identity but in general is smaller.

\begin{question}  What conditions on the Hamiltonian
 $H: M\to {\R}$ imply that 
$\rho(\phi_t^H) \to \infty$ as $t\to \infty$?
\end{question}

Clearly, we need to assume that $H$ does not generate a circle action, or 
more generally, that its flow $\phi_t^H, t\ge 0,$ 
is not quasiperiodic in the sense 
that its elements  are not contained in any compact subset of $\Ham(M)$.
However, this condition is not sufficient.  For example, if $H$ has the 
form $F - F\circ\tau$ where $F\in \LieHam(M)$ has support in a set $U$ 
that is disjoined by $\tau$ (i.e. $\tau(U)\cap U = \emptyset$), then  
$$
\phi_t^H \;=\; [\phi_t^F, \tau]  \;=\;  \phi_t^F \tau (\phi_t^F)^{-1} \tau^{-1}.
$$
Hence 
$$
\rho(\phi_t^H) \;\le\; \rho(\phi_t^F \tau (\phi_t^F)^{-1}) + 
\rho(\tau^{-1})\; \le\; 2\rho(\tau).
$$
Therefore  $\rho(\phi_t^H)$ is bounded.  But it is easy to construct 
examples on $S^2$ for which the sequence $\phi_n^H, n = 
1,2,3,\dots$ has no subsequence that converges in the $C^0$ topology.
Thus the flow is not quasiperiodic.
\MS

\begin{exercise}\rm  Suppose that $M$ is a Riemann surface of genus $>0$ 
and that $H$ has a level set that represents a nontrivial element in 
$\pi_1(M)$.  Show that  $\rho(\phi_t^H) \to \infty$  as $t\to \infty$ by 
lifting to the universal cover and using the energy-capacity inequality.
\end{exercise}

\begin{question} Is the sum $\rho^+ + \rho^-$ of the one sided seminorms
nondegenerate for all $(M, \om)$?
\end{question}

The seminorms $\rho^+$ and $\rho^-$ are said to be \lq\lq one sided" 
because they do not in general take the same values 
on an element $\phi$ and its 
inverse.  (As we show in Proposition~\ref{prop:polt}
they also have very natural geometric 
interpretations in which they measure the size of $\phi$ from 
just one side.)  Their sum $\rho^+ + \rho^-$ is two sided.  Hence its null 
set
$$
\nnull(\rho^+ + \rho^-) = \{\phi:\rho^+(\phi) + \rho^-(\phi) = 0\}
$$
is a normal subgroup of $\Ham(M)$.  Therefore it is trivial: in other words,
$\rho^+ + \rho^-$ is either identically zero or is nondegenerate.
Thus to prove nondegeneracy one just has to find one element on which
$\rho^+ + \rho^-$ does not vanish.
The paper~\cite{McvH} develops geometric arguments 
(using Gromov--Witten invariants on suitable 
Hamiltonian fibrations over $S^2$) 
that  show that it is nondegenerate in certain cases,
for example if $(M, \om)$ is a projective space or is weakly exact, 
i.e. $\om$ vanishes on $\pi_2(M)$.  (For the latter case, see also 
Schwarz~\cite{Sch}.) 
However the general case is still open.  An even harder question is 
whether the one sided norms $\rho^{\pm}$ are each nondegenerate.  Now 
the null set is only a conjugation invariant semigroup  and so could 
be a proper subset of $\Ham(M)$.  Hence it seems that one could only prove
nondegeneracy by  an argument that would apply to an arbitrary element of
$\Ham(M)$.

\subsection{Questions of stability}\label{ss:stab}

One of the  first nonelementary results in symplectic topology is due to 
Eliashberg~\cite{ELWAVE} and Ekeland--Hofer~\cite{EKELHOF}: 
 \MS
 
 \NI
 $\bullet$ {\bf Symplectic rigidity theorem:}
{\it the group $\Symp(M,\om)$ is $C^0$-closed in $\Diff(M)$.}\MS
 
 This celebrated result  
 is the basis 
 of  symplectic topology. The proof shows that there are
 invariants $c(U)$ (usually called {\bf symplectic capacities}) of an open 
 subset of a symplectic manifold that are
continuous with respect to the Hausdorff metric on sets and that are
preserved only by symplectomorphisms. (When $n$ is even,
one must slightly modify the previous statement to rule out 
the case $\phi^*(\om) = -\om$.)  There are 
several ways to define  suitable $c$.   Perhaps the easiest
is to take Gromov's width:
$$
c(U) = \sup \{\pi r^2: B^{2n}(r)\mbox{ 
embeds symplectically in }U\}.
$$
Here  $B^{2n}(r)$ is the standard ball of radius $r$ in Euclidean 
space.\MS

\NI
 $\bullet$   {\it Stability properties of $\Symp(M)$ and 
 $\Symp_0(M)$.}\MS
 
 By this we mean that if $G$ denotes either of these
 groups, there is a $C^1$-neighbourhood $\Nn(G)$ of $G$
 in $\Diff(M)$ that deformation retracts onto $G$.  This follows from 
 the Moser isotopy argument. In the case $G = \Symp(M)$ take
 $$
 \Nn(\Symp)= \{\phi\in \Diff(M): t\phi^*(\om) + (1-t)\om \mbox{ is 
 nondegenerate}, for all t\in [0,1]\}.
 $$
 By Moser, one can define for each such $\phi$ a unique  isotopy 
 $\phi_t$ (that 
 depends smoothly on $\phi^*(\om)$) such that 
 $\phi_t^*(t\phi^*(\om) + 
 (1-t)\om) = \om$ for all $t$.  Hence $\phi_1\in \Symp(M)$.
 Similarly, when $G= \Symp_0(M)$ one can take $ \Nn(G)$ to be the
 identity component of $\Nn(\Symp)$.
  \MS

\begin{question}  Is there a $C^0$-neighborhood of $\Symp_0(M)$ in 
$\Diff_0(M)$ that deformation retracts into $\Symp_0(M)$?
\end{question}

\NI
   $\bullet$ {\it The Flux conjecture}\MS
   
The analogous questions for the Hamiltonian group are much harder, 
and much less is known.  We do not even know if the group $\Ham(M)$ 
is always closed in the $C^1$-topology, let alone whether it is 
$C^0$-closed.  
It follows from diagram~(\ref{eq:flux}) that $\Ham(M)$ is $C^1$-closed 
in $\Symp_0(M)$ if and only if the flux group  $\Ga_\om$ is closed.
Therefore, the most important question here is the following.

\begin{question}  Is $\Ga_\om$ a discrete subgroup of $H^1(M, \R)$?
\end{question}  

The hypothesis that it is always discrete is known as the {\bf Flux conjecture.}
$\Ga_\om$ is known to be discrete for many $(M,\om)$.
For example, it follows from~(\ref{eq:flux0}) that $\Ga_\om$ is a 
subgroup of $H^1(M; {\mathcal P}_{\om})$ where ${\mathcal P}_{\om}$ 
denotes the periods of $[\om]$, i.e. its values on 
$H_2(M;\Z)$.\footnote
{
In fact, one can restrict to its set of values on spheres by 
Proposition~\ref{prop:lm0} (i) below.}
Thus $\Ga_\om$ is discrete 
if  $[\om]\in H^2(M;\Z)$.  Another easy case is when
$$
\wedge [\om]^{n-1}: H^1(M,\R)\to  H^{2n-1}(M,\R)
$$
is an isomorphism (eg. if $(M,\om)$ is K\"ahler.) 
Nevertheless, this question
does not yet have a complete answer: 
see~Lalonde--McDuff--Polterovich~\cite{LMP1,LMP2} and 
Kedra~\cite{Ked1,Ked2}.  It follows that there may 
be manifolds $(M, \om)$ for which the normal subgroup $\Ham(M)$ of 
$\Symp(M)$ is not closed with respect to the $C^\infty$ topology.  In 
fact, if $\Ga_\om$ is not discrete then one should think of $\Ham(M)$ as a 
leaf in a foliation of $\Symp_0(M)$ that has codimension equal to the first 
Betti number $\rk(H^1(M,\R))$ of $M$.\MS

\NI
$\bullet$  {\it Special geometric properties of elements in 
   $\Ham(M)$.}\MS
   
One indication that $\Ga_\om$ may always be discrete is that the 
elements $\Ham(M)$ have special geometric properties.
   In particular, according to  Arnold's celebrated conjecture 
   (proven by  Fukaya--Ono and Liu--Tian based on work by Floer)
  the number of fixed points of $\phi\in\Ham$ may be estimated as
  $$
  \#{\rm Fix\,} \phi\ge \sum_k {\rm rank\,} H^k(M, \Q)
  $$
  provided that all its fixed points are nondegenerate (i.e. the graph  
  of $\phi$ is transverse to the diagonal.)
The following simple argument shows how this is related to the Flux 
problem.  Denote  $r(M): \sum_k {\rm rank\,} H^k(M, \Q).$

\begin{lemma}  Suppose that $\La\subset H^1(M;\R)$ has the property 
that every  element $\al\in H^1(M;\R)\setminus \La$ is the 
flux of some symplectic path $\{\phi_t^\al\}_{t\in [0,1]}$ whose 
time-$1$ map is nondegenerate and has $< r(M)$ fixed points.  
Then $\Ga_{\om} \subset \La$.
\end{lemma}
\NI
{\bf Proof:} Suppose this is false, and let $\{\psi_t\}$ be a 
loop in $\Ham(M)$ with flux equal to $\al\in H^1(M;\R)\setminus \La$.   
Then the path $\{\psi_t^{-1}\phi_t^\al\}_{t\in [0,1]}$  is 
Hamiltonian and has time $1$-map $\phi_1^\al$.  Therefore $\phi_1^\al$
must have at least $r(M)$ fixed points by Arnold's conjecture, 
which contradicts the 
hypothesis.
\QED\MS

For example, one can apply this to the standard  torus $(T^{2n}, \om_0) = 
(\R^{2n}/\Z^{2n},\om_0)$, taking $\La$ 
to be the lattice $H^1(T^{2n};\Z)$, and conclude that the flux group 
must be contained in this lattice.  But  we knew this anyway. 
So far, no one has succeeded in getting very far with this kind of 
geometric argument.\MS

\NI
$\bullet$ {\it Stability of Hamiltonian loops}\MS

Although what one might call geometric stability (even for fixed $\om$)
has not yet been 
established for the Hamiltonian group, it does have stability 
properties on the homotopy level.  Here is a typical question.  
Suppose given  continuous map from a finite CW complex $X$ to $\Ham(M,\om)$.
What happens if we perturb $\om$?   
  If $\om'$ is sufficiently close to $\om$, then it follows by Moser 
  stability that $X$ is homotopic through maps to $\Diff(M)$  to a map 
  into $\Symp(M, \om')$.  But can we always deform $X$ into $\Ham(M, 
  \om')?$
Since $\pi_k(\Ham(M)) = \pi_k(\Symp_0(M))$ when $k>1$ by
diagram~(\ref{eq:flux}) this is automatic when $X$ is simply connected.
However the case $X = S^1$ is not at all obvious, and is proved 
in~\cite{LMP2,Mcq}.  In the version stated below $(M^M)_{id}$ denotes the group of 
homotopy self-equivalences of $M$, i.e. the identity component 
of the space of degree $1$ maps $M\to M$.

\begin{prop} Suppose that $\phi\in \pi_1(Symp(M, \om))$ and 
$\phi'\in \pi_1(Symp(M, \om'))$ represent the same element of 
$\pi_1((M^M)_{id})$.  Then 
$$
\Flux_{\om}(\phi) = 0 \quad\IFF\quad \Flux_{\om'}(\phi') = 0.
$$
\end{prop}

This is an easy consequence of a more general vanishing theorem for 
various actions of the Hamiltonian  group.  One can generalise the 
defining equation
(\ref{eq:flux0}) for the Flux homomorphism to get a general 
definition of  the action $\tr_\phi: H_*(M)\to H_{*+k}(M)$
of an element $\phi\in H_k(M^M)$.  Namely  
 if $\phi$ is represented by the cycle $t\mapsto \phi_t$ for $t\in V^k$
and $c\in H_*(M)$ is represented by $x\mapsto c(x)$ for $x\in C$
then $\tr_\phi(c)$ is represented by the cycle
$$
V^k\times Z \to M: (t,x)\mapsto \phi_t(x).
$$
To say this action is trivial means that
$$
\tr_\phi(c) = 0\qquad\mbox{whenever } c\in H_k(M), \;k > 0.
$$
If $V^k=S^k$ is a sphere then this action is precisely the 
differential $\p$ in 
the Wang long exact sequence for the associated bundle $M\to P\to
S^{k+1}$:
$$
\dots\to H_i(M)\to H_i(P) \to  H_{i-k}M\stackrel{\p}\to H_{i-1}M\to\dots.
$$
Therefore the triviality of the action is equivalent to saying that 
$H_*(P)$ is isomorphic to the tensor product $H_*(M)\otimes H_*(S^{k+1})$.  

The results to date on these questions are still incomplete.
The following proposition states the most important known conclusions.
The first result below is a consequence of
the proof of the Arnold conjecture.  Another more direct proof may be found in
Lalonde--McDuff--Polterovich~\cite{LMP2}.  The second part is proved 
in Lalonde--McDuff~\cite{LMh} and the third is an easy consequence.
Since the main ideas in the proofs are sketched in the survey 
article McDuff~\cite{Mcox} we shall not say more about them here.

\begin{prop}\label{prop:lm0}
\begin{description}\item[(i)]
The evaluation map $\pi_1(\Ham(M)) \to \pi_1(M)$ is zero.
\item[(ii)] The natural action of $H_*(\Ham(M),\Q)$ on $H_*(M,\Q)$ is 
trivial.
\item[(iii)]  If $(M, \om)\to P\to S^{k+1}$ is any bundle with 
structural group $\Ham(M)$ then $H_*(P;\Q)\cong H_*(M;\Q)\otimes 
H_*(S^{k+1};\Q)$.
\end{description}
\end{prop}

It would be very interesting to know if there are any similar 
results in other related categories, such as the category of Poisson 
manifolds.  Even though it does not seem as if there would be a good 
analogue for the Hamiltonian group as such, there still might be 
a notion of something akin to a Hamiltonian bundle and hence some 
analogue of property (iii) above.  

Note that (iii) definitely fails 
for general symplectic bundles over $S^2$: if $\{\psi_t\}_{t\in S^1}$ is a 
symplectic loop with nontrivial flux then the corresponding bundle 
over $S^2$ has nontrivial Wang differential.  This situation is 
discussed further in Section~\ref{sec:lect4}.  But when $k>1$ there is 
no difference between Hamltonian and symplectic bundles over $S^{k+1}$, 
since, as we already mentioned, $\pi_k(\Ham) = \pi_k(\Symp)$ for 
$k>1$.

Here is another rather curious consequence of Proposition~\ref{prop:lm0}.
It applies to {\it any} symplectic fibration $(M, \om)\to P\to S^2$, 
not just Hamiltonian fibrations. 

\begin{cor}  Let $(M, \om)\to P\to S^2$ be any symplectic fibration 
and let $\p$ be its Wang differential.  Then $\p\circ\p: H_iM\to 
H_{i+2}(M)$ is zero.
\end{cor}

It is quite possible that this is always true for any smooth 
fibration over $S^2$ (with compact total space), but I do not know a 
proof. 

\section{The homotopy type of $\Symp(M)$}

\subsection{Ruled $4$-manifolds}\label{ss:rr} 

Another set of questions concerns the homotopy type of the group 
$\Symp(M)$.  
In rare cases this is completely understood.  The following 
results are due to Gromov~\cite{GRO}:

\begin{prop}\label{prop:gro}
\begin{description}\item[(i)]
 $\Symp^c(\R^4, \om_0)$ is contractible;
\item[(ii)]
$\Symp(S^2\times S^2, \si + \si)$ is homotopy equivalent to the extension 
of $\SO(3)\times \SO(3)$ by $\Z/2\Z$ where this acts by interchanging the 
factors;
\item[(iii)]  $\Symp(\C P^2, \om)$ is 
homotopy equivalent to $\PU(3)$.
\end{description}
\end{prop}

It is no coincidence that these results occur in dimension $4$.  The 
proofs use $J$-holomorphic spheres, and these give much more information in 
dimension $4$ because of positivity of intersections.  
Abreu~\cite{Abr} and Abreu--McDuff~\cite{AM} recently extended Gromov's 
arguments to other ruled surfaces.  
Here are their main results, stated for convenience for the product 
manifold $\Sigma\times S^2$ (though there are similar results for the 
nontrivial $S^2$ bundle over $\Si$.)   We sketch the easiest
proofs in Section~\ref{sec:lect3}. 

Consider the following family of symplectic forms on $M_g = \Sigma_g\times S^2$
(where $g$ is  genus$(\Si)$):
$$
\om_\mu = \mu \si_\Si +  \si_{S^2},\qquad \mu > 0,
$$
where $\si_Y$ denotes (the pullback to the product of) an area form on
the Riemann surface $Y$ with total area $1$.\footnote
{
Using results of Taubes and Li--Liu, Lalonde--McDuff show in~\cite{LMc} that 
these are the {\it only} symplectic forms on $\Si\times S^2$ up to 
diffeomorphism.}
Denote by $G_\mu = G_\mu^g$ the subgroup  
$$
 G_\mu^g:= \Symp(M_g,\om_\mu)\cap \Diff_0(M_g)
$$
of the group of symplectomorphisms of $(M_g,
\om_\mu)$.  When $g > 0$ $\mu$ ranges over all positive numbers. 
However, when $g=0$ there is an extra symmetry ---  interchanging the two
spheres gives an isomorphism $G_\mu^0 \cong G_{1/\mu}^0$ --- and so we 
take $\mu \ge 1$.  Although it is not completely obvious, 
there is a natural homotopy class of maps  from  $G_\mu^g$ to $
G_{\mu + \eps}^g$ for all $\eps>0$.  To see this, let
$$
G_{[a,b]}^g = \bigcup_{\mu\in [a,b]} \; \{\mu\}\times G_\mu^g 
\;\;\subset \;\;\R\times \Diff(M_g).
$$
It is shown in~\cite{AM} that the inclusion $G_b^g\to G_{[a,b]}^g$
is a homotopy equivalence. Therefore we can take the map
 $G_\mu^g \to G_{\mu + \eps}^g$ to be the composite of the inclusion
$G_\mu^g\to G_{[\mu, \mu+\eps]}^g$ with a homotopy inverse 
$G_{[\mu, \mu+\eps]}^g\to G_{\mu+\eps}^g$.

\begin{prop}\label{prop:am1}  As $\mu\to \infty$,  the groups $G_\mu^g$
tend to a limit $G_\infty^g$ that
has the homotopy type  of the identity component $\Dd_0^g$ of the group 
of fiberwise diffeomorphisms
of $M_g = \Si_g\times S^2 \to \Si$.
\end{prop}

\begin{prop}\label{prop:am2}
When $\ell  < \mu \le \ell + 1$ for some integer
$\ell \ge 1$,  
$$
H^*(G_\mu^0, \Q) = \La(t, x, y)\otimes \Q[w_\ell], 
$$
where $\La(t,x,y)$ is an exterior algebra over $\Q$ with generators $t$ of
degree $1$, and $x,y$ of degree $3$ and $\Q[w_\ell]$ is the polynomial
algebra
on a generator $w_\ell$ of degree $4\ell$.
\end{prop}

In the above statement, the generators $x,y$ come from 
$H^*(G_1^0) = H^*(\SO(3)\times
\SO(3))$ and $t$ corresponds to an element in $\pi_1(G_\mu^0), \mu >
1$, found by Gromov in~\cite{GRO}.   
Thus the subalgebra $\La(t, x, y)$ is the pullback of $H^*(\Dd_0^0, \Q)$
under the map $G_\mu \to \Dd_0^0$.    The other generator
$w_\ell$ is fragile, in the sense that the corresponding element in 
homology disappears (i.e. becomes null
homologous) when $\mu$ increases. 

There are still many unanswered questions about these groups.  Here is a 
sampling.

\begin{question}
Is the group $G_\mu^g$ connected for all $\mu > 0, g\ge 0$?
\end{question}

It is shown in McDuff~\cite{Mcrs} that the answer is \lq\lq yes" whenever $\mu \ge 
[g/2]$.  (The case $g=0$ was proved in~\cite{AM}.)  
This paper also 
provides an affirmative answer to
 the next question in the genus zero case.

\begin{question}
Is the homotopy type of the groups $G_\mu^g$ constant on the intervals
$\mu\in (\ell, \ell+1]?$ 
\end{question}

\begin{question}
The group $G_\mu^1$ is known to be constant for $0< \mu \le 1.$  What is 
its homotopy type?
\end{question}

The methods used to prove the above results extend to certain other 
closely related manifolds.  For example, Pinsonnault in   
thesis~\cite{Pin} studies 
the symplectomorphism group of the one point blow up 
of $(S^2\times S^2, \om_\mu)$ with $\mu=1$ and shows it is homotopic 
to the $2$-torus, $T^2$.  As shown by Lalonde--Pinsonnault in~\cite{LALPi}
this group becomes more complicated when $1< \mu\le 2$, and its 
homotopy groups change when the blow up radius $r$ passes through the 
critical level $\pi r^2 = \mu - 1.$  This implies that the homotopy 
type of the associated 
space of symplectically embedded balls also changes, the first known 
example of such a phenomenon.

There has  been some attempt to generalize these results to higher 
dimensions. Le--Ono~\cite{LO} and 
Buse~\cite{Bu} use parametric 
Gromov--Witten invariants to
obtain information on the symplectomorphism groups 
of products $(M_g, \om_\mu)\times (N,\om)$, 
while Seidel studies the case of 
products of two projective spaces in~\cite{SEI}. 
\MS

\subsection{The topology of $\Symp(M)$ for general $M$}

\NI
$\pi_0(\Symp(M))$

This group is known as {\bf  the symplectic mapping class 
group}.  Seidel has done interesting work here, studying
{\it symplectic Dehn twists} especially on manifolds with 
boundary.  He considers  the group $\Symp(M, \p M)$ of 
symplectomorphisms that are the identity near the boundary, 
detecting quite large subgroups of $\pi_0(\Symp(M, \p M))$
by using Floer homology
to study the effect of Dehn twists on the
Lagrangian submanifolds in $M$: cf~\cite{Seid2}.
\MS

\NI
$\pi_1(\Symp(M))$

This is an abelian group and one can try to detect its elements by studying 
various natural homomorphisms.  One such is the Flux homomorphism:
$$
\Flux_\om: \pi_1(\Symp(M)) \to \Ga_\om \;\;\subset \;\;H^1(M,\R),
$$
that has kernel equal to $\pi_1(\Ham(M))$.
There are several other interesting homomorphisms defined on this kernel,
most notably a homomorphism to the units in the quantum homology ring of 
$M$ known as the Seidel representation~\cite{Seid,LMP1}: cf. \S4.
\MS

Very little is known about the higher homotopy groups.  Observe, however,
 that the existence of the diagram~(\ref{eq:flux})
implies that the inclusion $\Ham(M) \to \Symp(M)$ induces an isomorphism 
on $\pi_j, j>1$.  \MS

\subsection{Characteristic classes}

Reznikov shows in~\cite{Rez} how to define classes $\la_k\in 
H^{2k-1}(\Ham(M), \R)$ for $k > 1$ that he calls higher Cartan classes
 using the invariant inner product of~(\ref{eq:inner})
and an analog of Chern--Weil theory.  
However, one can define the corresponding classes 
$$
c_k^H\in H^{2k}(B\Ham(M),\R), k\ge 2,
$$
on the classifying space using the notion of the coupling class of a 
symplectic fibration.  As we will see in~\S3, given any smooth fibration
$\pi:M\to P\to B$ with structural group $\Ham(M,\om)$ there is a canonical 
class $u\in H^2(P)$ that extends the symplectic class on the fibers.
$u$ is called the Guillemin--Lerman--Sternberg {\bf coupling class} 
and is characterized by the property that 
the integral of $u^{n+1}$ over the fibers of $\pi$ is $0$ in $H^2(B)$:
cf~\cite{GLS,Pbk}, or~\cite[Ch~6]{MS1}.  Then we set
$$
c_k^H: = \int_M u^{k+n},\quad k>1.
$$
This defines $c_k^H\in H^{2k}(B)$.  By naturality this has to come from a 
class $c_k^H\in H^{2k}(B\Ham(M))$ that we shall call
a {\it Hamiltonian Chern class}.

Other characteristic classes  can be constructed using  
the Chern classes of the vertical tangent bundle  $T_{vert}P\to P$ whose 
fiber at $x\in P$ is the tangent space to the fiber through $x$.  Denoting 
these classes by $c_k^{vert}$ we get corresponding elements in 
$H^*(B\Ham(M,\om))$ by integrating products of the form
$$
u^\ell c_{k_1}^{vert}\dots c_{k_p}^{vert}
$$
over the fibers of $P\to B$. Rather little is known about these 
classes, though
Januszkiewicz and Kedra~\cite{JK} have recently 
calculated them for symplectic toric manifolds.
(They appear in slightly different guise in~\cite{LMP2}.  Note also that 
the classes $c_k^{vert}$ exist for symplectic bundles, and so when 
$\ell=0$ the classes extend to $H^*(B\Symp(M))$.)

Now consider a situation in which a compact Lie group acts on $(M,\om)$ 
in such a way as to induce an injective homomorphism $G \to \Ham(M)$.
(Such an action is often called weakly Hamiltonian.)  Then there is a 
corresponding map of classifying spaces:
$$
BG \to B\Ham(M)
$$
and one can ask what happens to the Hamiltonian Chern classes $c_k^H$ under 
pullback.  It follows from Reznikov's definitions that if $G$ is 
semisimple and the action is effective then the pullback of $c_2^H$ is nonzero.  Hence 
$$
H^3(\Ham(M))\ne 0
$$
whenever $(M,\om)$ admits such an action.  Moreover
he shows by direct calculation
that, in the case of the action of the projective unitary 
group $\PU(n+1)$ on complex projective space $\C P^n$, the pullbacks of the
$c_k^H$ to $B\PU(n+1)$ are multiplicatively independent.  
Thus the inclusion $\PU(n+1)\to \Ham(\C P^n)$ 
induces an injection on rational homotopy.

\begin{question}  What conditions imply that
such an inclusion $G \to \Ham(M)$ is nontrivial homotopically?
\end{question}

The above question is deliberately vague:  what does one mean precisely by 
\lq\lq nontrivial"?  Presumably one could extend Reznikov's 
calculation to
other actions of compact semisimple simply connected Lie groups $G$ on 
their homogeneous K\"ahler quotients: see Entov~\cite{En}.  

Reznikov's argument is elementary.  In contrast, the next result
uses fairly sophisticated analytic 
tools: cf. McDuff--Slimowitz~\cite{MSlim}.  

\begin{prop}  Given a semifree Hamiltonian action of $S^1$ 
on $(M,\om)$, the associated homomorphism
$$
\pi_1(S^1) \to \pi_1(\Ham(M))
$$
is nonzero.
\end{prop}

\NI
{\it Proof.}
Recall that an action $\{\phi_t\}_{t\in \R/\Z}$ of $S^1$ is called semifree if 
no stabilizer subgroup is proper.  Equivalently, the only points fixed by 
some $\phi_T$ for $0 < T < 1$ are fixed for all $t$.  Therefore, by the 
remarks after Question~\ref{q:slim}, both the paths $\{\phi_t\}_{t\in [0,3/4]}$
and $\{\phi_{-t}\}_{0\le t\le 1/4}$ are length minimizing in their 
respective homotopy classes.  Since they have different lengths, these 
paths cannot be homotopic. 
\QED

 This result has 
been considerably extended by McDuff--Tolman~\cite{MT}.  Note also 
that the semifree hypothesis is crucial: for example the action
of $S^1$ on $M=\C P^2$ given by
$$
[z_0:z_1:z_2]\mapsto [e^{2\pi i\theta} z_0: e^{-2\pi i\theta}z_1:z_2]
$$
gives rise to a nullhomotopic loop in $\Ham(M)$. On the other hand the 
image of $\pi_1(S^1)$ in $\pi_1(\Ham(M))$  might be finite: for example
the  rotation of $S^2$ by one turn has order $2$ in $\Symp(S^2)\simeq 
\SO(3)$.

In an ongoing 
project~\cite{KMc},  Kedra and McDuff have recently 
shown that  if a Hamiltonian circle 
action  is inessential (i.e. gives rise to a contractible element 
in $\pi_1(\Ham(M))$) then there is an associated nontrivial element 
in $\pi_3(\Ham(M))$.    This extends Reznikov's result: 
in the case when the circle is a subgroup in 
a semisimple Lie group $\G$ then this element is precisely the one  he
detected. 
\MS

Finally, in~\cite{LO} Le--Ono define {\bf Gromov--Witten 
characteristic classes} on $B\Symp(M)$.   If one restricts to
$B\Ham(M)$ and fixes the genus $g$, these are 
indexed by the elements of $H_2(M;\Z)$
and can be defined as follows.  Let $M_{\Ham}\to B\Ham(M)$
denote the bundle with fiber $M$ associated to the obvious action of 
$\Ham(M)$ on $M$.  Then, by Proposition~\ref{prop:lm0}
 $H_2(M_{\Ham})$ splits as the sum $H_2(M) \oplus 
H_2(B\Ham(M))$.  Hence 
each $A\in H_2(M;\R)$ gives rise to a well defined class in
the homology of the fibers of the fibration $M_{\Ham}\to B\Ham(M)$, and
 there is a class $\GW_A\in 
H^{\mu(A)}(B\Ham(M))$ whose value on a cycle $f:B\to B\Ham(M)$
is the \lq\lq number of isolated $\TJ$-holomorphic genus $g$ curves in 
class $A$" in the pullback 
fiber bundle $f^*(M_{\Ham}) \to B$.  Here the almost complex structure $\TJ$
is compatible with the fibration in the sense that it restricts on each 
fiber $(M_b,\om_b)$ to a tame almost complex structure, and the index
$$
\mu(A) = (g-1)(2n-6) - 2c_1(A).
$$  
(Of course, to define the invariants correctly one has to 
regularize the moduli space in the usual way: see for 
example~\cite{RuT}.)

\begin{question}  {\it When are
these characteristic classes nontrivial?}
\end{question}
 
Unfortunately one cannot get very interesting examples from the groups 
$G_\mu^0$ discussed above.
Le--Ono show that when $\mu>1$ the $1$-dimensional vector space 
 $H^2(BG_\mu^0, \R)$ 
is generated by a Gromov--Witten 
characteristic class $\GW_A$.  However, the nontriviality of this 
class $\GW_A$  can also be proved by purely 
homotopical means since the corresponding loop
in $G_\mu^0$ does not vanish in the group of 
self homotopy equivalences of $S^2\times S^2$.  (In fact,
Le--Ono show that the cohomology ring
of the corresponding fibration over $S^2$ is not a product ring.)
It is shown in~\cite{AM} that the new elements $w_\ell\in H^*(G_\mu^0)$
do not transgress to  $H^*(BG_\mu^0)$, but rather give rise to 
{\it relations} in this ring.  

One can also define classes by evaluating 
appropriate moduli spaces of $J$-holomorphic curves at $k$ points for 
$k > 0$.  One interesting fact pointed out in Kedra~\cite{Ked2} is 
that the nontrivial Gromov--Witten classes of dimension 
$\mu(A) = 0$ constrain the image of the Flux homomorphism.

\subsection{$J$-holomorphic curves in $S^2\times S^2$}\label{sec:lect3}

In this section we shall give a brief overview of the proof of some 
of the results on $G_{\mu}: = \Symp(M\times M,\om_\mu)$ mentioned in 
\S\ref{ss:rr}.  Fuller details may be found in the survey article~\cite{LM1}
as well as in Abreu's beautiful paper~\cite{Abr}.
 An introduction to some of the 
technicalities may be found in the lecture notes~\cite{Mcut} or
 the more broadly based survey~\cite{Mcell}.  Proposition~\ref{prop:gro} 
is proved in exhaustive detail in~\cite{MS2}.
The proofs are based on the 
behavior of $J$-holomorphic spheres in $4$-manifolds.  Here is a 
summary of their most important properties.\MS

\subsection{Analytic background}

\NI
{\bf Almost complex structures:}\,  An almost complex structure $J$ on 
$M$ is an automorphism of $TM$ with square $-1$.  It is said to be 
$\om$-tame if $\om(v,Jv) > 0$ for every nonzero  $v\in TM$.  For every 
symplectic manifold, the space $\Jj(\om)$ of $\om$-tame almost complex 
structures is nonempty and contractible.  Given $J\in \Jj(\om)$ we 
shall always use the associated Riemannian metric 
$$
 g_J(v,w): =\frac 12 (\om(v,Jw) + \om(w,Jv)
 $$
on $M$.
 \MS

\NI
{\bf $J$-holomorphic curves:}\,  A map $u:(S^2,j)\to (M,J)$ is said to be 
$J$-holomorphic if it satisfies the following nonlinear 
Cauchy--Riemann equation:
$$
\pbar_J u: = \frac 12 (du + J\circ du\circ i) = 0.
$$
If $J\in \Jj(\om)$ then $\om$ restricts to an area form at all non 
singular points in the image of $u$.  Therefore if $u$ represents the 
homology class $A\in H_2(M)$ we must have $\om(A) > 0$ unless $A=0$ 
and $u$ is constant. Note also that the M\"obius group $\PSL(2,\C)$ 
acts by reparametrization $u\mapsto u\circ\ga$ 
on the space of $J$-holomorphic spheres.  

In many respects the 
behaviour of these curves is exactly the same as in the integrable case.
In particular, in dimension $4$ there is {\bf positivity of 
intersections}: every intersection point of two distinct $J$-holomorphic curves
contributes positively to their intersection number, with
nontransverse intersections contributing $> 1$.  Thus if $A\cdot B=0$ 
every $J$-\hol $A$-curve is disjoint from every  $J$-\hol $B$-curve, 
while if $A\cdot B=1$ 
every $J$-\hol $A$-curve meets every  $J$-\hol $B$-curve precisely 
once and transversally.  This is relatively easy to prove when one is 
considering two distinct curves.  A rather more subtle result is that 
this remains true for a single (non multiply covered) curve; in 
particular any singular point on a curve $u$ (i.e. point where $du=0$)
contributes positively to the self-intersection number. This can 
be formulated as the {\bf adjunction inequality}:
\MS

\NI
$\bullet$ {\it if $u:(S^2,j)\to (M^4,J)$ is a $J$-\hol sphere in class $A$ 
then $j$ is an embedding if and only if }
$$
c_1(A) = 2 + A\cdot A.
$$
Thus there is a homological criterion for a curve to be embedded. In 
particular if $A$ can be represented by a $J$-\hol sphere then {\it every} 
$J'$-\hol representative of $A$ is embedded for any $J'\in \Jj(\om)$.
 \MS

\NI
{\bf The moduli space $\Mm(A;J)$:}\, The main 
object of interest to us is the moduli 
space $\Mm(A;J)$ of all $J$-holomorphic $A$-spheres.  This has two 
fundamental  properties.
\MS

\NI
$\bullet$  {\it There is a set $\Jreg$ of second category in $\Jj(\om)$ 
such that $\Mm(A;J)$ is a smooth manifold of dimension $2n+2c_1(A)$ 
for every $J\in \Jreg$.}

The proof involves Fredholm theory.
\MS

\NI
$\bullet$  Since $\PSL(2,\C)$ is noncompact, the space $\Mm(A;J)$ is 
never compact except in the trivial case when $A=0$.  However the 
quotient $\Mm(A;J)/\PSL(2,\C)$ is sometimes compact. Even if it is 
not, it has a well behaved compactification that is made up of genus 
$0$ stable maps in class $A$.   These objects were called 
cusp-curves by Gromov and are sometimes known as bubble trees.  One 
shows that $\Mm(A;J)/\PSL(2,\C)$ can be noncompact only if $A$ has a 
representation by a connected union of two or more $J$-holomorphic spheres.
In particular, this can happen only if the class $A$ decomposes as a
a sum $A= A_1+A_2$ where $\om(A_1)> 0, \om(A_2)> 0$.  Hence
$\Mm(A;J)/\PSL(2,\C)$ is compact in the case when $[\om]$ is integral 
and $\om(A) = 1$.
\MS

\NI
{\bf Fredholm theory and regularity:\,}   
We now say a little more about Fredholm theory and regularity since 
this is so crucial to our argument.\MS

Let $\Mm(A,\Jj)$ be the space of all pairs $(u,J)$, where $u:(S^2,j)\to (M,J)$
is $J$-\hol,
$u_*([S^2]) = A\in H_2(M)$, and $J\in \Jj(\om)$. One shows that a suitable 
completion of $\Mm(A,\Jj)$ is a Banach manifold and that the projection 
$$
\pi: \Mm(A,\Jj)\to \Jj
$$
is Fredholm.  
In this situation one can apply an infinite dimensional 
version of Sard's theorem (due to Smale) that states that there is a set 
$\Jj_{reg}$ of second category in $\Jj$ consisting of regular values of 
$\pi$, i.e. points where $d\pi$ is surjective.  
Moreover by the implicit function theorem for Banach manifolds the inverse
image  of a regular value is a smooth manifold of dimension equal to the index of 
the Fredholm operator $\pi$. Thus one finds that for almost every $J$ 
$$
\pi^{-1}(J) = \Mm(A;J)
$$
is a smooth manifold of dimension $2(c_1(A) + n)$.  The index 
calculation here follows by investigating the linearization of $\pi$ 
which turns out to be essentially the same as the linearization
\begin{equation}\label{eq:Du}
D_u: C^\infty(S^2,u^*(TM))\to \Om_J^{0,1}(S^2,u^*(TM))
\end{equation}
of the operator $\pbar_Ju$.  One can check that $D_u$ is a zeroth 
order perturbation of the usual Dolbeault differential $\pbar$ from functions 
to $(0,1)$-forms with values in the bundle $u^*(TM)$. 
Hence $D_u$ has the same index as $\pbar$, which in turn  is
given by the Riemann--Roch theorem.  

There is another important point here.  When $u:S^2\to M^4$ is an 
embedding, the bundle $u^*(TM)$ splits (as a complex bundle) into the sum 
of the tangent bundle to $S^2$ with the normal bundle $L$ to $\im 
u$.  The restriction of $D_u$ to the tangent bundle is always 
surjective.  It follows that $D_u$ is surjective if and only if its 
restriction to the line bundle $L$ is surjective.  In this case, the
Riemann--Roch theorem says there is a dichotomy:\MS

$\bullet$
{\it either $c_1(L) \ge -1$ and $D_u$ is surjective;} 

$\bullet$
or {\it   $c_1(L) < -1$ and the rank of $\coker\, D_u$ is constant (and 
equal to $2(|c_1(L)| - 1)$.)}\MS

We will use this fact later.

Next observe that, by a transversality theorem for paths, given any two elements 
$J_0, J_1\in \Jj_{reg}$ there is a path $J_t, 0\le t\le 1$, such that the 
union 
$$
W = \cup_t\Mm(A;J_t) = \pi^{-1}(\cup_t J_t)
$$ is a smooth (and also 
oriented) manifold with boundary 
$$
\p W = \Mm(A;J_1)\cup - \Mm(A;J_0).
$$
  It follows 
that the evaluation map 
$$
ev_J: \Mm(A;J)\times_{PSL(2,\C)} S^2 \to M,\quad (u,z)\mapsto 
u(z),
$$
is independent of the choice of (regular) $J$ up to
oriented bordism.  In
particular, {\em if} we could  ensure that everything is compact and {\em if} we 
arrange that $ev$ maps between  manifolds of the same dimension then the
degree of this map would be  independent of $J$.
\MS

\subsection{The case $S^2\times S^2$}

In the case of $(S^2\times S^2, \om_\mu)$ we are interested in looking at curves
 in the classes $A: = [S^2\times pt]$ and $B: = [pt\times S^2]$.  
Thus $c_1(A) =c_1(B) = 2$ and
$$
\dim(\Mm(A;J)\times_G S^2) =  2(c_1(A) + 2) - 6 + 2 = \dim(M).
$$
If $\mu = 1$ the 
above remarks about compactness imply that 
 the  moduli spaces $\Mm(A;J)/PSL(2,\C))$ are always
compact  and $ev_J$ has degree $1$.  Moreover, for each $J$
the $J$-\hol $A$-curves are mutually disjoint (by 
positivity of intersections), and one can show that they form the 
fibers of a fibration of $S^2\times S^2$.  Similar statements holds 
for the $B$-curves.

When $\mu> 1$ this remains true 
for the smaller sphere $B$ (though one needs some extra arguments to 
prove this).   On the other hand, it is now possible for the $A$-curve 
to decompose since $\om_{\mu}(A-B) > 0$.  Moreover the class $A-B$ 
is represented by the symplectic embedding of $S^2$ onto
the antidiagonal  $z\mapsto (z,-z)$ where we think of $z\in S^2\subset\R^3\}$,
and this submanifold can be made $J$-holomorphic for suitable 
$J$.\footnote
{
Given a symplectically embedded submanifold $C$ of $(M, \om)$ one can always 
choose $J\in \Jj(\om)$ so as to preserve $TC$ since the set of 
choices at each point is contractible.  When $C$ is $2$-dimensional,
every almost complex structure on it is integrable.  Hence there is a 
complex structure $j$ on $C$ for which the obvious inclusion 
$\io:(C,J)\to (M, J)$ is $J$-\hol.}
Thus there are $J\in \Jj(\om_{\mu})$ for which the curve $A-B$ is 
represented.  Since $(A-B)\cdot(A-B) = -2$, positivity of 
intersections implies that  this representative is 
unique.  Moreover, there cannot be any $A$-curves since 
$A\cdot(A-B) = -1$.  It follows that when $1<\mu\le 2$ 
every $J\in \Jj(\om_\mu)$ is of one of two kinds:\MS

\NI
(i)  if $A$ is represented, $S^2\times S^2$ has two transverse
fiberings, one by $A$-curves and one by $B$-curves;

\NI
(ii) if $A$ is not represented, $S^2\times S^2$ is still 
fibered by the $B$-curves; however there is a unique $(A-B)$-curve 
and there are no $A$-curves.
\MS

Thus we may decompose $\Jj: = \Jj_{\om_\mu}$ into the disjoint union 
of two sets: $\Jj_0$ on which $A$ is represented and
$\Jj_1$ on which $A-B$ is represented.   (For $\mu\in (\ell-1, \ell]$
one defines  $\Jj_{k}$ for $k\le \ell$ to be the set of $J$ for 
which $A-kB$ is represented.)

\MS

To go further, we need to use the consequences of the Riemann-Roch theorem
that were mentioned above.  Since the normal bundle to an $A$- or 
$B$-curve is trivial, every $A$- and every $B$- curve is regular, 
i.e. $D_u$ is always onto in this case.  Since regularity is an open 
condition this means that the set $\Jj_0$ is open.  On the other hand 
the class $A-B$ cannot be represented by a regular curve $u$, 
i.e. $D_u$ can never be surjective, since $c_1(L) = -2$ in this case.
But the cokernel of $D_u$ has constant rank, and one can show that
this implies that $\Jj_1$ is a submanifold of $\Jj$ of codimension 
$2$.  (In fact the normal bundle to $\Jj_1$ in $\Jj_0$ at a point 
$(u,J)$ can be 
identified with $\coker D_u$: see~\cite{Abr}.)   Thus $\Jj(\om_\mu)$ is 
a stratified space when $1<\mu\le 2$ with $2$ strata. The picture is 
similar when $\mu> 2$ except that there are now more strata: see~\cite{Macs}. 
This 
complicates the calculations that we present below.   
\MS

\NI
{\bf The strata $\Jj_i$ as homogeneous spaces:\,} In what follows we 
either suppose that $\mu = 1$ and let $i = 0$ or 
suppose that $1<\mu \le 2$ and let $i = 0, 1$.

Each stratum $\Jj_i$ 
contains an integrable element $J_i$.  We may take $J_0$ to be 
the product structure $j\times j$ and $J_1$ to be the Hirzebruch 
structure obtained by identifying $S^2\times S^2$ with the 
projectivization  $\PP(L_2\oplus \C)$, where $L_2\to \C P^1$ is a 
holomorphic line bundle of Chern class $2$.  (Note that the section 
given by $\PP(L_2\oplus 0)$ is rigid, with self-intersection number 
$-2$: hence it corresonds to the antidiagonal.)   Denote by $K_i$ the 
maximal compact subgroup in the identity component of the
complex automorphism group of $J_i$.  Thus 
$K_0 \cong \SO(3)\times SO(3)$, while $K_1 \cong \SO(3)\times S^1$.  
Here the $\SO(3)$ factor can be identified with the diagonal subgroup 
of  $K_0$, while the $S^1$ factor is generated by an $S^1$ action.  
This can either be thought of as coming from the action $[w_0:w_1]\mapsto
[e^{2\pi i t}w_0:w_1]$, or can be explicitly described by the formula:
$$
\phi_t: S^2\times S^2\to S^2\times S^2, \quad (z,w)\mapsto 
(z,R_{z,t}w),
$$
where $R_{z,t}$ rotates the sphere about the axis through the points 
$\pm z$ by the angle $e^{2\pi i t}$.  Note that the diagonal and 
antidiagonal are fixed by each $\phi_t$.

Now consider the map 
$$
G_{\mu}:= \Ham(S^2\times S^2,\om_{\mu}) \to \Jj_i,\quad 
\phi\mapsto \phi_*(J_i).
$$
Since the stabilizer of $J_i$ is precisely $K_i$ this induces a 
quotient map
\begin{equation}\label{eq:qot}
q_i: G/K_i\to \Jj_i.
\end{equation}
The claim is that these maps $q_i$ are homotopy equivalences.  When $\mu = 
1$ so that $\Jj_0$ is contractible, this will follow if we produce a  
map  $s_0: \Jj_0\to G/K_0$ such that $s_0\circ q_0 \sim {\rm id}$.   The 
general argument is a little more complicated and may be found 
in~\cite{Abr}.    

To construct $s_0$ we proceed as follows.  Fix a point $x_0$ in 
$S^2\times S^2$ and for each $J\in \Jj_0$ 
let $C_A, C_B$ be the unique $A$- and $B$-curves through $x_0$.
We shall think of these as \lq\lq coordinate axes'' and of the families 
of $A$ and $B$ curves as corresponding to their parallel translates 
$w= const$ and $z=const$.   More precisely,
choose parametrizations $u: S^2\to C_A, v: S^2
\to C_B$, and denote by $C_{w,A}$ the unique $A$-curve through 
the point $v(w)\in C_B$ and  by $C_{z,B}$ the unique $B$-curve through 
$u(z)\in C_B$.
Then define the map
$\phi_J\in \Diff(S^2\times S^2)$ by setting
$$
\phi_J(y) = (z,w)\in S^2\times S^2,\quad\mbox{where } y = C_{z,B}\cap 
C_{w,A}.
$$
Because the fibrations by $A$ and $B$ curves are transverse,
this map $\phi_J$ is a diffeomorphism.  Though it is not quite 
symplectic, it turns out that it lies in the \lq\lq Moser 
neighborhood'' of $\Symp$, in other words, it can be canonically 
isotoped into $G_\mu$. If $\phi_J'$ denotes the endpoint of this 
isotopy, we define 
$$
s_0(J): = \phi_J'\, K_0 \in G_1/K_0.
$$
To check that $s_0$ is well defined one must investigate the effect 
of changing the parametrizations $u,v$.  Initially these are defined 
modulo $\PSL_*$, the subgroup of $\PSL(2,\C)$ consisting of elements that 
fix a point.  Since this is homotopy equivalent to $S^1$, it is not 
hard to see that there is a consistent choice of $u,v$ modulo 
$S^1\subset \SO(3)$, at least over compact families.  There are many 
ways of getting round this point, for example by 
using  balanced maps as in~\cite{MS2}.
\MS

\NI
{\bf Proof of Proposition~\ref{prop:gro}(ii).}\, When $\mu = 1$ there 
is just one stratum: $\Jj=\Jj_0$.  Since $\Jj$ is contractible, the 
result is immediate 
from~(\ref{eq:qot}).  
\QED\MS

\NI
{\bf Proof of Proposition~\ref{prop:am2} in the case $1<\mu\le 2$}.
In this case there are two strata so that $\Jj=\Jj_0\cup \Jj_1.$
By Alexander--Spanier duality we know that
$$
H^i(\Jj_1)\cong H^{i+1}(\Jj_0),\qquad i\ge 0.
$$
Abreu also shows that the spectral sequence of the fibration
$
K_i\to G_\mu\to \Jj_i$ degenerates for $i=0,1$.  Therefore 
$$\begin{array}{ll}
 H^*(G_{\mu};\R) \cong H^*(\Jj_0)\otimes H^*(\SO(3)\times\SO(3)) 
 &\quad (i)\\
H^*(G_{\mu};\R) \cong H^*(\Jj_1)\otimes H^*(S^1\times \SO(3))&\quad 
(ii).
\end{array}
$$
The argument is now a simple algebraic computation.  Since $G_\mu$ is 
a Hopf space, its cohomology algebra is free.  It follows 
that $H^*(\Jj_i)$ is also a free algebra.  Let $t$ (resp. $x.y$) be the image in 
$H^*(G_{\mu};\R)$ of
the generator of $H^1(S^1)$ (resp. the generators of 
$H^3(\SO(3)\times \SO(3))$.)  Comparing the ranks of $H^3$ in (i) and 
(ii) we see that $H^3(\Jj_1)$ must have at least one generator, say 
$x_1$.  Therefore $H^4(\Jj_0)$ has a corresponding generator, say $w$.
One now proves that there can be no other new generator:  if there were,
let $k$ be the minimum dimension of such an element in $H^*(G_\mu)$ 
and use the isomorphism $H^k(\Jj_0)\cong H^{k-1}(\Jj_1)$ to show that 
there would also have to be a new generator in dimension $k-1$.  Hence  
$H^*(G_{\mu};\R)$ is generated by $t,x,y,w$.\QED\MS

The proof for the other cases $\mu> 2$ is similar, but both its aspects 
become more complicated because there are more strata $\Jj_k$.  
One must work harder to show that the $\Jj_k$ do form a 
stratification of $\Jj$ (see~\cite{Macs}), and the algebraic 
calculation is also considerably more elaborate.

\section{Symplectic geometry of fibrations over $S^2$}
\label{sec:lect4}

Many of the proofs of the propositions above rely on 
properties of Hamiltonian fibrations over $S^2$.  In this lecture we 
consider the geometric properties of such fibrations, relating them to 
the Hofer norm.  The  main ideas in this section come from Lalonde--McDuff 
and Polterovich.

\subsection{Generalities}

 Consider a smooth fibration $\pi: P  \to 
B$ with fiber $M$, where $B$ is either  $S^2$ or the $2$-disc $D$.
Here we consider $S^2$ to be the union $D_+ \cup D_-$ of two copies of 
$D$, with the  orientation  of $D_+$.  We denote the equator $ D_+\cap 
D_-$ by $\p$, oriented as the boundary of $D_+$, and choose some point
$*$ on $\p$ as the base point of $S^2$.  Similarly, $B = D$ is provided 
with a basepoint $*$ lying on $\p = \p D.$  In both cases, we assume that 
the fiber $M_*$ over $*$ has a chosen identification with $M$.

Since every smooth fibration over a disc can be trivialized, 
we can build any smooth fibration $P\to S^2$ by taking two product 
manifolds $D_\pm\times M$ and gluing them along the boundary $\p\times M$ by a 
based 
loop $\la = \{\la_t\}$ in $\Diff(M)$.  Thus
$$
P = (D_+\times M)\; \cup \;(D_-\times M)/\sim,\quad 
(e^{it},\la_t(x))_+\equiv (e^{it},x)_-.
$$
A symplectic fibration is built from a based loop in $\Symp(M)$ and a 
Hamiltonian fibration from one in $\Ham(M)$.  Thus the smooth fibration
$P\to S^2$ is symplectic if and only if there is a smooth family of 
cohomologous symplectic 
forms $\om_b$ on the fibers $M_b$.   
It is shown in~\cite{Seid,MS1,LMh} 
that a symplectic fibration $P\to S^2$ is Hamiltonian if and only if
the fiberwise forms $\om_b$ have a closed extension $\Om$.  
(Such forms $\Om$ are called {\bf $\om$-compatible}.)
 Note that, in any of these categories,
two fibrations are equivalent if and only if their defining loops are 
homotopic.

From now on, we restrict to Hamiltonian fibrations.
By adding the pullback of a 
suitable area form on the base we can choose the closed extension
 $\Om$ to be  symplectic. 
Observe that there is a unique class $u\in H^2(P, 
 \R)$  called the {\bf coupling class} 
that restricts to $[\om]$ on the fiber $M_*$ and has the property that
$$
\int_P u^{n+1} = 0.
$$
(This class will have the form $[\Om] - \pi^*(a)$ for a suitable $a\in 
H^2(B, \p B)$.)  Correspondingly we decompose $\Om$ as
\begin{equation}\label{eq:coup}
\Om = \tau + \pi^*(\al)
\end{equation}
where $[\tau] = u$, and call $\tau$ the {\bf coupling form}.
(Although we will not need this, there is a canonical choice for
the form  $\tau$ depending only on the connection defined by $\Om$.)

The closed form $\Om$ defines a
connection on $\pi$ whose horizontal distribution is $\Om$-orthogonal 
to the fibers.   
If $\ga$ is any path in $B$ then $\pi^{-1}(\ga)$ is 
a hypersurface in $P$ whose characteristic foliation consists of 
the horizontal lifts of $\ga$,
and it is not hard to check that
the resulting holonomy is Hamiltonian round every
contractible loop, and hence round every loop. 
(A proof is given in~\cite[Thm~6.21]{MS1}.)   

Thus, given $\Om$, $\pi$ can be symplectically 
trivialized over each disc $D_\pm$
 by parallel translation along 
a suitable set of  rays.  This means that there is a fiber preserving 
mapping  
$$
\Phi_\pm:  \pi^{-1}(D) \to M\times D_\pm, \quad \Phi|_{M_*} = id_M
$$
such that the pushforward $(\Phi_\pm)_*\Om$ 
restricts to the same form $\om$ on each fiber $M\times pt$.
These two trivializations differ by  a loop
$$
e^{it}\mapsto \phi_t = \Phi_+\circ (\Phi_-)^{-1}(e^{it})\in \Symp(M,\om)
$$
where $e^{it}$ is a coordinate round the equator $\p = D_-\cap D_+$.

\begin{exercise}\rm  Check that this loop is homotopic to
the defining loop for the fibration $P\to S^2$.
(Since $\pi_1(G)$ is abelian for any group $G$ it does not matter 
whether or not we restrict to based homotopies.)
\end{exercise}

\begin{defn} 
The {\bf monodromy}  $\phi = \phi(P)\in \Ham(M)$  
of a fibration $(P, \Om)\to B$
is defined to be the monodromy of the connection determined by 
$\Om$ around the based oriented loop $(\p, *)$.
 Using the trivialization of
$P$ over $\p$ provided by $B$ itself if $B = D$
or by $D_+$ if $B = S^2$, one gets  a well defined lift 
$\tphi$ of $\phi$ to the universal cover $\THam$.  
\end{defn}

\begin{exercise}\rm  Start with a fibration $(P,\Om)\to S^2$ and 
break it in half to get two fibrations $(P_+,\Om)\to D_+$
and  $(P_-,\Om)\to D_-$.  Here the inclusion $D_+\to S^2$ is orientation 
preserving, while $D_-\to S^2$ is orientation reversing.
What is the relation between the monodromies of these fibrations,
(a) considered as elements in $\Ham(M)$ and (b) considered as elements in 
$\THam(M)$?
\end{exercise}

In Section~\ref{ss:ABW} we shall consider a fibration over  a base 
$B$ which is the sphere
$S^2$ with some points removed.  We assume that 
near each deleted point  (i.e. end 
of $B$) the fibration is identified with the product $[0,\infty)\times 
S^1\times M$ and that the
form $\Om$ is normalized so that in the coordinates $(s,t,x)$ it can 
be written as
 $a(s,t) ds\wedge dt + \om - d_MH_t\wedge dt$, where $d_M$ 
denotes the exterior derivative on $M$.

\begin{exercise}\label{ex:Hend}\rm  Check that the monodromy of $\Om$ round such an 
end is precisely the Hamiltonian flow of $H_t: = H_{t+1}$. (The sign 
conventions are given in~(\ref{eq:H}).)
\end{exercise}

\subsection{The area of a fibration}

We define the {\it area} of a fibration $(P, \Om) \to B$ 
to be:
$$
{\rm area}\,(P, \Om): = \frac {{\rm vol\,}(P,\Om)}{{\rm\vol\,}(M, \om)}
= \frac{\int_{P}\Om^{n+1}}{(n+1)\int_{M}\om^{n}}. 
$$
Thus a product fibration $(B\times M, \al_B + \om)$ has area $\int_B 
\al_B$.

\begin{exercise}\label{ex:area}\rm
Decompose $\Om$ as $\tau + \pi^*(\al)$ as in~(\ref{eq:coup}).
Show that $\area(P,\Om) = \int_B \al$.
\end{exercise}

The next definition describes ways to use this area to measure
the size of elements of $\Ham(M)$ or $\THam(M)$.

\begin{defn}\label{def:area}
\begin{itemize}
\item[(i)]
 $\ta^+(\tphi)$ (resp. $a^+(\phi)$)
is the infimum of $\area\,(P, \Om)$ taken over 
all $\om$-compatible symplectic forms $\Om$ on the fibration $P\to D$
 with monodromy $\tphi$ (resp. $\phi$).
\item[(ii)] 
$a(\phi)$ is the infimum of $\area\,(P, \Om)$ taken over {\it all} 
fibrations $(P, \Om)\to S^2$ with monodromy $\phi$.
\item[(iii)]  $\ta^-(\tphi): = \ta^+(\tphi^{-1})$ and 
$a^-(\phi): = a^+(\phi^{-1})$.
\end{itemize}
\end{defn}

We now show that these area measurements
agree with Hofer type measurements.  Recall that
$\rho^+(\phi)$ is the infimum of $\Ll^+(H_t)$ 
over all Hamiltonians $H_t$ whose times $1$ map is $\phi$.
Similarly, we define  $\trho^+(\tphi)$ to be the 
infimum of $\Ll^+(H_t)$ over all Hamiltonians $H_t$ whose flow over 
$t\in [0,1]$ is a representative of the element $\tphi\in \THam(M)$.

The following lemma is a slightly sharper version of
Polterovich's results in~\cite{Peg}. 

\begin{prop}\label{prop:polt}\begin{itemize}
\item[(i)]  $\trho\,^{+}(\tphi) =
\ta^+(\tphi)$;
\item[(ii)]  $\rho^{+}(\phi)  + \rho^-(\phi) = a(\phi)$;
\end{itemize}  
\end{prop}

We prove (i).  The proof of (ii) follows easily (see~\cite{McvH}).

\MS

\NI
{\bf Proof that $\trho\,^{+}(\tphi) \ge \ta^+(\tphi)$:}\MS

This is by direct 
construction.  Suppose that the path $\phi_t^H$ generated by 
$H_t$ is $\tphi$.  For simplicity let us suppose that the functions 
$\min(t) = \min_x H_t(x)$ and $\max(t) = \max_x H_t(x)$  are
smooth, so that by replacing $H_t$ by 
$H_t - \max(t)$ we have that $\max_x H_t(x) = 0$ for all $t$.  
Suppose also that $H_t(x) = 0$ for all $x\in M$ and all $t$ 
sufficiently near $0,1$.\footnote
{
This can be arranged without altering the time-$1$ map or the 
Hofer length.  For this and other 
technicalities see~\cite{LM} or~\cite{McvH}.}
Then define
the graph $\Ga_H$ of $H_{t}$ by
$$
\Gamma_H: = \{ (x, t, H_t(x)): x\in M, t\in [0,1] \} \subset M \times
[0,1]  \times {\R}.
$$ 
For some small $\eps> 0$ choose a smooth function
$\mu(t): [0,1] \rightarrow [0,+2\eps]$ such that
$$
\int_0^1 \mu(t)\, dt = \eps.
$$
Consider the following {\bf thickening of the region over} $\Ga_H$: 
$$
R_H^+(\eps): = \{ (x,t,h) \; \vline \; H_t(x) \leq h \leq \mu(t)\} 
\subset M \times [0,1]  \times {\R}.
$$ 
Note that if $\mu$ is chosen for $t$ near  $0,1$ to be tangent to the 
lines $t = const$ at $t = 0,1$ we may
arrange that $R_{H}^{+}(\eps)$ is a manifold with corners along $t= 
0,1$.
(Recall that $H_{t}\equiv 0$
for $t$ near $0,1$.) 

Let $\Om_0 = \om + dt\wedge dh$ be the standard symplectic form on 
$M \times [0,1]  \times {\R}$.  Then the monodromy of the hypersurface 
$\Ga_H$ (oriented as the boundary
of $R_H^+(\eps)$) is precisely $\phi_t^H$, while the rest of the boundary 
has trivial monodromy.  Further, it is easy to define a projection
$\pi$ from $R_H^+(\eps)$ to the half disc  $HD$ whose fibers all lie in
the hypersurfaces $t=const$.  Thus, after rounding the corners, we get a 
fibered space $\pi: R_H^+(\eps) \to D$ with monodromy $\tphi$.
It remains to check that the area of  $(R_H^+(\eps), \Om_0)$ (before 
rounding corners)
is precisely $\Ll^+(H_t) + \eps$. \QED

\begin{rmk}\rm
Similarly, we can define a
manifold with corners $R_H^-(\eps)$ that thickens
 the region below $\Ga_H$ by setting 
$$
R_H^-(\eps): = \{ (x,t,h) \; \vline \; \min(t) - \mu(t) \leq h \leq H_t(x) \} 
\subset M \times [0,1]  \times {\R}.
$$
Note that $\area (R_H^-(\eps)) = \Ll^-(H_t)+\eps$.
\end{rmk}

To prove the other inequality we combine 
Polterovich's arguments from~\cite{Pk}\S3.3 and~\cite{Peg}\S3.3.

\begin{lemma}\label{le:polt3} $\ta^{+}(\tphi) \ge \trho\,^{+}(\tphi)$.
\end{lemma}
\proof{}  
Suppose we are given a fibration $(P,\Om)\to D$ 
with area $ < \trho\,^{+}(\tphi)$ and monodromy $\tphi$.
By Moser's theorem we may isotop $\Om$ so that it is a product 
in some neighborhood $\pi^{-1}(\Nn)$ of the
base fiber $M_*$. Identify the 
base $D$ with the unit  square $K = \{0\le x,y \le 1\}$ taking 
$\Nn$ to a neighborhood of $\p' K = \p K - \{1\}\times (0,1)$, 
and then identify $P$ with 
$K\times M$ by parallel translating along the lines $y=$ const.
% $\{(x,y): x\in [0,1]\}$. 
In these coordinates, the form $\Om$ may be written as 
$$
\Om\;\; = \;\; \om +  d_MF'\wedge dy + L' dx\wedge dy
$$
where $F', L'$ are suitable functions on $K\times M$ and $d_M$ denotes 
the fiberwise exterior derivative.  
Because $\Om$ is a product near  $\pi^{-1}(\p' K)$, 
$d_MF' = 0$ there and $L'$ reduces to a function of $x,y$ only. By subtracting 
a suitable function $c(x,y)$ from $F'$  we can arrange that $F = F' - c(x,y)$ 
has zero mean
 on each fiber $\pi^{-1}(x,y)$ and then write 
$L' + \p_xc(x,y)$  as $-L + 
 a(x,y)$ where 
$L$ also has zero fiberwise means.  Thus 
\begin{equation}\label{eq:form}
\Om\;\; = \;\; \om +  d_MF\wedge dy - L dx\wedge dy + a(x,y)dx\wedge dy,
\end{equation}
where both $F$ and $L$ vanish near $\pi^{-1}(\p' K)$ and have zero  
fiberwise means.
 Since $\Om$ is symplectic it must be nondegenerate on the
 $2$-dimensional distribution \,\,{\it Hor} formed by the 
 $\Om$-orthogonals to the fibers.  Hence we must have 
 $-L(x,y,z)+ a(x,y) > 0
$ for all $x,y\in K, z\in M$.  Moreover, because $L$ has zero 
fiberwise means, 
$\area (P, \Om) = \int a(x,y) dx \wedge dy.$  Hence
\begin{equation}\label{eq:a2}
\int \max_{z\in M} L(x,y,z) \,dx\wedge dy < \int a(x,y)\,dx\wedge dy = 
\area(P,\Om).
\end{equation}

We claim that $-L$ is the curvature of the induced
connection $\Om_\Ga$.  To see this, consider the vector fields
$X = \p_x, Y = \p_y - \sgrad F$ on $P$ that are
 the  horizontal lifts of $\p_x, \p_y$.\footnote
{
Here the symplectic gradient  $\sgrad F$ is defined by setting
  $\om(\sgrad F, \cdot) = -d_MF(\cdot)$.}
It is easy to check that their commutator $[X,Y]= XY - YX$ is 
vertical and that
$$
[X,Y] = -\sgrad(\p_xF) = \sgrad L 
$$
on each fiber $\pi^{-1}(x,y)$ as claimed.  (In fact, the first three terms 
in~(\ref{eq:form})  make up the coupling form $\tau_{\Ga}$.)

Now let $f_{s}\in \Ham(M)$ be the monodomy of $\Om_\Ga$ along the path
$t\mapsto (s,t), t\in [0,1]$.  (This is well defined because all  
fibers have a natural identification with $M$.)
The path $s\mapsto f_{s}$ is a Hamiltonian isotopy
from the identity to $\phi = f_{1}$, and it is easy to see that it is 
homotopic to the original path $\tphi$ 
given by parallel transport along $t\mapsto 
[1,t]$.  (An intermediate path $p_T$ in the homotopy 
might consist of the 
path $s\mapsto f_{s}^{T}$ for $0\le s\le T$, where $f_{s}^{T}$ is 
the monodromy along $t\mapsto (s/T,t), t\in [0,T]$, followed by the lift
of $\tphi_s, s\in [T,1]$, to $\THam$.)  Therefore 
 $\Ll^+(f_s) \ge \trho^+(\tphi)$, and we will derive a contradiction 
 by estimating $\Ll^+(f_s) $.

To this end, let $X^s, Y^t$ be the (partially defined) flows of the vector 
fields $X,Y$ on $P$ and set $h_{s,t} = Y^t X^s$.  Consider the 
$2$-parameter family of (partially defined)
vector fields $v_{s,t}$ on $P$ where $v_{s_0,t_0}(p)$ is tangent to the path 
$s\mapsto h_{s,t_0}$ at $h_{s_0,y_0}(p)$ for $p = (x,y,z)\in P$.  Thus
$$
v_{s,t} = \p_s h_{s,t} = Y^t_*(X)\quad\mbox{on }\;\; {\rm Im\,}h_{s,t}.
$$
In particular $v_{s,1}(x,y,z)$ is defined when $y=1, s\le x$.
Since $f_s = h_{s,1}$ we are interested in calculating the 
vertical part of $v_{s,1}(s,1,z)$.  Since the points 
with $y = 1$ are in ${\rm Im\,}h_{s,t}$ for all $(s,t)$ we may write
\begin{eqnarray*}
v_{s,1} & = &  v_{s,0}+ \int_0^1 \p_t(v_{s,t})\; dt \\
& = &  \p_x + \int_0^1 Y^t_*([X,Y])\; dt .
\end{eqnarray*}
We saw above that $[X,Y] = \sgrad L$. Hence $Y^t_*([X,Y]) = \sgrad(L\circ 
(Y^t)^{-1})$ and
\begin{eqnarray*}
v_{s,1}(s,1,z) & = & \p_x + \int_0^1 \sgrad(L((Y^t)^{-1}(s,1,z)) \; dt\\
& = & \p_x + \sgrad \int_0^1\,L(s,1-t, (Y_v^t)^{-1}(z))\;dt 
\end{eqnarray*}
where $Y_v^t$ denotes the vertical part of $Y^t$.
Hence the  Hamiltonian $H_{s}$
that  generates the path $f_{s}, s\in [0,1],$ and has zero mean
satisfies the inequality
$$
H_{s} (z)\; \le\; \int \left(\max_{z\in  M} L(s,t,z)\right) dt \;\le 
\; \int a(s,t)  dt
$$
since $L(s,t,z) \le a(s,t)$ by~(\ref{eq:a2}).
Thus $\trho\,^+(\tphi)\le \area P$, contrary to hypothesis.\QED

\section{The quantum homology of fibrations over $S^2$}\label{sec:lect5}

In this lecture we show how to use the Seidel representation
$$
\Ss: \pi_1(\Ham(M,\om))\to (\QH_{ev}(M))^{\times}
$$
to estimate areas of fibrations and also get information
on the homotopy properties of Hamiltonian fibrations.
\MS

\subsection{Quantum and Floer homology}

First of all, what is the small quantum homology ring $\QH_*(M)$?
Because it is more efficient, we shall use the formulation 
in~\cite{MS2}, which is slightly different from~\cite{LMP2,Mcq}.

Denote by 
$\Lambda^\univ$ the universal Novikov ring formed by
all formal power series with rational coefficients $\la_\kappa$ of the form
$$
     \lambda = \sum_{\kappa\in\R}\la_\kappa t^\ka,\qquad
     \#\left\{\kappa\in\R\,|\,\la_\kappa\ne0,\,
     \kappa\ge c\right\} < \infty\mbox{ for all }c\in\R.
$$
Thus the power $\kappa$ of $t$ is allowed to go to $-\infty$. 
Set $\La: = \Lambda^{\univ}[q,q^{-1}]$, where $q$ is a variable of 
degree $2$.
% 
% 
% Set $c_1 = c_1(TM)\in H^2(M,\Z)$.  Let $\La$ be the   Novikov 
% ring of
% the group $\Hh = H_2^S(M,\R)/\!\!\sim$  with valuation $I_\om$ where 
% $B\sim B'$
% if $\om(B-B') = c_1(B-B') = 0$.
%  Thus $\La$ is the completion of the rational group ring
%  of $\Hh$ with elements of the form 
% $$
% \sum_{B\in \Hh} q_B\; e^B
% $$
% where for each $\ka$ there are only finitely many nonzero
% $q_B\in \Q$ with $\om(B) > - \ka$.
Additively, the quantum homology  is simply the usual homology with 
coefficients in $\La$:
$$
\QH_{*}(M) := \QH_*(M;\La)
= H_*(M)\otimes\La.
$$
We may define an $\R$ grading on $\QH_*(M;\La)$ by setting
$$
\deg(\al\otimes q^dt^{\kappa}) = \deg(\al) + 2d,\qquad \al\in H_j(M),
$$
but can also think of
$\QH_*(M;\La)$ as $\Z/2\Z$-graded with 
$$
\QH_{\ev} =  
H_{\ev}(M)\otimes\La, \quad \QH_{odd} =  
H_{odd}(M)\otimes\La.
$$

The quantum intersection product  is linear over $\La$ and is defined 
on classes $\al\in H_i(M), \be\in 
H_j(M)$ as follows.  We abbreviate $c: = c_1(M)$.
% $$
% a*b\in \QH_{i+j - 2n}(M), \qquad \mbox {for }\; a\in H_i(M), b\in 
% H_j(M)
% $$
%  is linear over $\La$ and is defined on the generators defined as follows:
\begin{equation}\label{eq:qm0}
     \al* \be = \sum_{B\in H_2^S(M)} (\al*\be)_B\otimes 
     q^{-c(B)}\,t^{-\om(B)}\;\in \;\QH_{i+j - 2n}(M),
\end{equation} 
where  $(\al*\be)_B\in H_{i+j- 2n+2c_1(B)}(M)$ is defined by the 
requirement
that 
\begin{equation}\label{eq:qm}
(\al*\be)_B\,\cdot\, \ga = \GW^M_{B,3}(\al,\be,\ga) \quad\mbox{ for all 
}\;\ga\in 
H_*(M).
\end{equation}
Here  $\GW^M_{B,3}(\al,\be,\ga)$ denotes the Gromov--Witten invariant that counts 
the number of $B$-spheres in $M$ meeting representing cycles for the
classes $\al,\be,\ga\in 
H_{*}(M)$,
and we have written $\cdot $ for the usual  intersection
pairing on $H_*(M) = H_*(M, \Q)$. In good cases, one can compute 
$\GW^M_{B,3}(\al,\be,\ga)$ as the intersection number of the class 
$\al\times \be\times \ga$ in $M^3$ with the image of the evaluation map
$$
\ev:\Mm(M,B;J)\to M^3,\quad u\mapsto \bigl(u(0),u(1), u(\infty)\bigr),
$$
where $J$ is a generic $\om$-tame almost complex structure and
$\Mm(M,B;J)$ denotes the $(2n+2c(B))$-dimensional moduli space of $J$-holomorphic 
$B$-spheres in $M$ as discussed in Section~\ref{sec:lect3}.
(In general one needs to use the virtual moduli cycle.)
Note that 
$\al\cdot  \be = 0$ unless $\dim(\al) + \dim (\be) = 2n$ in which case it is 
the algebraic number of intersection points of the cycles.
The product  $*$ is
extended to $\QH_*(M)$ by linearity over $\La$, and is associative. 
Moreover, it  preserves the $\R$-grading.  

 This product $*$ gives $\QH_{*}(M)$
the structure of a
graded commutative ring with unit $\1 = [M]$. Further, the invertible
elements in $\QH_{\ev}(M)$ form a commutative group  
$(\QH_{\ev}(M,\La))^\times$
that acts on   $\QH_*(M)$ by  quantum multiplication. 
\MS

We shall say very little about Floer homology; basic definitions can 
be found for example  in~\cite{Mcell,SAL3,MS0,MS2}.  
It is the Morse complex of the action 
functional on the loop space of $M$.  Given a Hamiltonian function 
$H_h: = H_{t+1}$ satisfying a suitable nondegeneracy hypothesis,
one forms a chain complex that is generated as a 
$\La$-module by the $1$-periodic orbits $\bx$ of $H$.  When $H$ is 
autonomous (i.e. independent of time) and is sufficiently $C^2$-small, 
these orbits are simply the critical points of $H$.  The 
$\by$-coefficient of the boundary 
 $\delta(\bx) $ 
is formed by counting isolated {\bf Floer trajectories} from $\bx$ 
to $\by$, 
i.e. solutions $u:\R\times S^1\to M$ of the {\bf Floer equation}
\begin{equation}\label{eq:Fl}
\p_su + J(u)(\p_t(u) - X_H(u)) = 0, \quad
\lim_{s\to -\infty}u(s,\cdot) = \bx,\;
\lim_{s\to \infty}u(s,\cdot) = \by.
\end{equation}
The resulting homology groups are denoted $\FH_*(M;H,J)$.  They are 
all canonically isomorphic.
For each (nondegenerate)  $H$, there is  a canonical isomorphism $\Phi$ from quantum homology  
$\QH_*(M;\La)$ to Floer homology $\FH_*(M;H,J)$ given by counting 
maps $u:\C\to M$ that are $J$-holomorphic near the unit disc $D$
with Floer boundary conditions at infinity; in other words, if we 
identify $\C\setminus D$ with the cylinder $(0,\infty)\times S^1$ 
then $u(s,t)$ satisfies the Floer equation~(\ref{eq:Fl}) for $s> 2$.
These isomorphisms are called PSS isomorphisms after 
Piunikhin--Salamon--Schwarz~\cite{PSS}.  This construction is 
important in the discussion of the ABW inequalities.

\subsection{The Seidel representation $\Ss$}

Now consider the fibration $P_\la\to S^2$ constructed from a loop $\la\in 
\pi_1(\Ham(M))$
as in \S\ref{sec:lect3}.
The manifold $P_\la$ carries two canonical 
cohomology classes,
the first Chern class of the vertical
tangent bundle 
$$ 
c: = c_1(TP_\la^{vert})  \in H^2(P_{\la},\Z),
$$
and the   coupling class  $u_\la$, i.e.  the unique class
in $H^2(P_\la,\R)$ such that 
$$
i^*(u_\la) = [\om],\qquad u_\la^{n+1} = 0,  
$$
where $i: M\to P_\la$ is the inclusion of a fiber.
We denote by $H_2^\si(P)$ set of {\bf section classes} in 
$H_2(P_\la,\R)$, i.e. the classes that project to the positive 
generator of $H_2(S^2;\Z)$.
% The next step is to choose a canonical (generalised) section class.
% By  section class, we mean a class that projects onto the positive generator
% of $H_{2}(S^{2}, \Z)$.  
% We consider the coset $\Aa = \si+H_2^S(M;\R)$ of all such classes,
% with  
% equivalence relation induced from that on $H_2^S(M;\R)$.
%  In the general case, 
% when $c_{1}$ and $[\om]$ induce linearly 
% independent homomorphisms $H_2^S(M) \to \R$, $\si_\la\in \Aa$ is
% defined by the requirement that
% \begin{equation}\label{eq:seccl}
% c_{vert}(\si_{\la}) \;= u_{\la}(\si_{\la})\; =\;0,
% \end{equation}
% which has a unique solution modulo the given equivalence.
% We show in~\cite{Mcq} that when $M$ is weakly exact
% such a class $\si_\la$ still exists and moreover is integral.
% In the remaining case, when $c_1$ is some multiple of $[\om]$ in
% $H_2^S(M)$,
% we choose $\si_\la$ so that $c_{vert}(\si_\la) = 0$. 
We then define the {\bf Seidel element}
\begin{equation}\label{eq:qm2}
\Ss(\la): = \sum_{\TB\in H_2^\si(P)} \al_\TB\otimes q^{-c(\TB)}t^{-u_\la(\TB)}
\end{equation}
where, for all $\ga\in H_{*}(M)$,
\begin{equation}\label{eq:qm3}
\al_{\TB}\cdot_{M} \ga = \GW^{P_{\la}}_{\TB,3}([M], [M], \ga).
\end{equation}
Note that $\Ss(\la)$ belongs to the strictly commutative part 
$\QH_{ev}$ of $\QH_{*}(M)$.
Moreover $\deg(\Ss(\la)) = 2n$.
It is shown in~\cite{Mcq} (using ideas from Seidel~\cite{Seid}
and Lalonde--McDuff--Polterovich\cite{LMP2}) that
for all $\la_1,\la_2\in \pi_1(\Ham(M))$
$$
\Ss(\la_1\la_2) = \Ss(\la_1)*\Ss(\la_2),\qquad 
\Ss(0) = \1,
$$
where $0$ denotes the constant loop.  Therefore 
$\Ss(\la)$ is invertible for all $\la$ and we get a 
representation
$$
 \Ss: \pi_{1}(\Ham(M,\om))\;\;\to\;\; (\QH_{\ev}(M;\La))^{\times}.
$$
Moreover since all $\om$-compatible forms are deformation 
equivalent, $\Ss$ is independent of the choice of $\Om$.
It is often useful to identify $(\QH_{\ev}(M;\La))^{\times}$ with the 
space $\Aut(\Q_*(M;\La))$ of automorphisms of $\QH_*(M;\La)$ 
considered as a (left) module over itself via the correspondence
$\al \mapsto \al*\cdot$.  We denote by $\Psi$ the 
resulting representation:
$$
\Psi: \pi_{1}(\Ham(M,\om))\;\;\to\;\; \Aut(\QH_*(M;\La)),\qquad
\Psi(\la)(a) : = \Ss(\la)*a.
$$
In this language, $\Ss(\la) = \Psi(\la)(\1)$.
\MS\MS

The fact that $\Ss(\la)$ is a unit means that Hamiltonian 
fibrations over $S^2$ always have plenty of holomorphic sections.   This has 
consequences for the  cohomology of $P_{\la}$. Indeed  one deduces that
 the map $H_*(M;\R)\to H_*(P_\la;\R)$ is injective by showing that
 any class $\be$ in its kernel must be annihilated by $\Phi(\la)$, 
 i.e. $\Phi(\la)(\be)=0$. Since $\Phi(\la)$ is an isomorphism, this implies $\be=0$.
 Since these arguments are sketched in~\cite{Mcox}, we shall 
 concentrate here on explaining other applications.

\subsection{Using $\Ss$ to estimate area}

Now let $\Om$ be any $\om$-compatible symplectic form on $P_\la$. 
As in Exercise~\ref{ex:area},
 its cohomology class has the form
$$
[\Om] = u_{\la} + \pi^{*}([\al])
$$ 
where  $\area(P_{\la}, \Om)  =  \int_{S^2} \al$.  The next results are due 
to Seidel.

Consider the valuation $v:\QH_{*}(M) \to \R$ defined by
\begin{equation}\label{eq:val}
v(\sum \al_{d,\kappa}\otimes q^dt^{\kappa}) = \sup \{\kappa: \al_{d,\kappa}\ne 
0\}).
\end{equation}
It follows from the definition of the quantum intersection product
in~(\ref{eq:qm0}), (\ref{eq:qm}) that
$v(\al*\be) \le v(\al) + v(\be).$  In fact, the following stronger statement is 
true. 

Denote the usual intersection product by  $\cap$, so that $a*b - 
a\cap b$ is the quantum correction to the usual product. 
Define
\begin{equation}\label{eq:Hbar}
\hbar = \hbar(M) = \min\left(\{\om(B): B\ne 0, \mbox{ some 
}\;\GW^M_{B,3}(\al,\be,\ga)\ne 
0\}\right),
\end{equation}
and note that $\hbar(M)> 0$:  standard compactness results 
imply that for each $c>0$ there are only finitely many classes $B$ with
 $\om(B) \le c$ that can be represented by a $J$-holomorphic curve 
 for generic $J$, and it is only such classes that give rise to 
 nonzero invariants.  If all 
the invariants $\GW^M_{B,3}(\al,\be,\ga)$ with $\al,\be,\ga\in H_*(M)$ and 
$B\ne 0$ vanish, we set $\hbar = \infty$.

\begin{lemma}\label{le:hbar}  For all $\al,\be\in
\QH_*(M)$,  $
v(\al*\be - \al\cap \be) \le  v(\al) +v(\be) - \hbar(M)$. 
\end{lemma}
\proof{}  This follows immediately from the definitions.\QED

\begin{prop}[Seidel]\label{prop:main} For each loop $\la$ in $\Ham(M)$  
$$
 \area(P_{\la}, \Om) \;\;>\;\; v(\Ss(\la)).
$$
\end{prop}
\proof{} 
% 
% Since $[\om]$ and $c_{1}$ are linearly independent on 
% $H_2^{S}(M)$ we may define $\Psi$ using a section $\si_{\la}$ that 
% satisfies~(\ref{eq:seccl}).
Let $\Ss(\la) = \sum_{d,\kappa} \al_{d,\kappa}\otimes q^dt^\ka,$
and let $v(\Ss(\la)) = \ka_0$.  Then $\al_{d,\kappa_0}\ne 0$.
By definition $\al_{d,\kappa_0}$ is determined by 
a count of $J$-holomorphic curves in 
$(P_{\la}, \Om)$ in the class $\TB$ where $c(\TB) = -d, u_\la(\TB) = 
-\ka_0.$  Hence 
this moduli space cannot be empty. Therefore
$$
\begin{array}{lcl}
0 \;\;< \;\;[\Om](\TB)  & = &  \pi^{*}([\al])(\TB) + u_\la(\TB)
\\
& = &  \int_{S^{2}}\al - \kappa_0 \\
&=& \area(P_{\la}, \Om) - \kappa_0.
\end{array}
$$
Thus $\area(P_{\la}, \Om) > \kappa_0$ for all $\Om$, as claimed.\QED

\begin{cor}\label{cor:main} In these circumstances 
$\trho^+(\la) \ge v(\Ss(\la))$.
\end{cor}
\proof{}  Combine
Proposition~\ref{prop:polt}(i) with Proposition~\ref{prop:main}.\QED

For applications of this estimate, see~\cite{McvH} and 
Lemma~\ref{le:ns1} below.

\subsection{The ABW inequalities}\label{ss:ABW}

We now give a very brief sketch of Entov's explanation of the ABW 
(Agnihotri--Belkale--Woodward) 
inequalities concerning the eigenvalues of products of unitary 
matrices.   Let us first consider the simplest case $G = \SU(2)$.
The question is: suppose we know the eigenvalues of the elements
$A_1,A_2\in G$.  What can we say about the eigenvalues of $A_1A_2$?
Equivalently, suppose we know that $A_1A_2A_3 = \1$.  What can we say 
about their eigenvalues?  In this case, $A_j$ has eigenvalues
 $e^{\pm 
2\pi i \zeta_j}$ for a unique $\zeta_j\in [0,1/2]$, and one can prove 
by computation that the triples $(\zeta_1,\zeta_2,\zeta_3)$ 
corresponding to solutions of the identity $A_1A_2A_3 = \1$ form the 
convex polygon in $[0,1/2]\times[0,1/2]\times[0,1/2]\subset \R^3$ 
described by the inequalities
$$
\zeta_1+\zeta_2+\zeta_3\le 1,\quad \zeta_1\le\zeta_2+\zeta_3,
\quad \zeta_2\le\zeta_3+\zeta_1, \quad \zeta_3\le\zeta_1+\zeta_2.
$$
There are corresponding inequalities that describe the relations 
among the  eigenvalues of matrices $A_j\in \SU(n), j=1,\dots,N$ such 
that $\prod_jA_j=\1$.
  Belkale found a generating set of 
inequalities for the resulting convex set.   
 Agnihotri--Woodward then observed that one may choose
 generators that are in bijective correspondence to the 
 nontrivial correlators in the quantum Schubert calculus.  (These 
 are the correlators that generate the relations in the quantum 
 cohomology of the Grassmannian.)
 We now give 
 a brief sketch of Entov's explanation of this fact.

 Entov considers the obvious action 
of $\SU(n)$ on the Grassmannian $M = \Gr(r,n)$ of complex $r$ planes 
in $\C^n$.  This  action is Hamiltonian.  Moreover each $A\in 
G\subset\Ham(M)$ is the 
time-$1$ map of a Hamiltonian flow generated by a function $H_A: 
M\to \R$ whose critical points  $x_I$ correspond 
to the fixed points of the action of $A$ on $\Gr(r,n)$ and hence to 
the $r$-planes that are spanned by eigenvectors.  Thus the label $I$
is a subset  of $\{1,\dots,n\}$  of 
cardinality $r$.    Moreover the
critical value of $H_A$ at $x_I$ is just
$\sum_{k\in I}\zeta_{Ak}$,  where $\zeta_{A1}\ge\dots\ge\zeta_{An}$ are 
the eigenvalues of $A$.
These critical values are precisely the terms that occur in the 
general ABW inequalities. (When $n = 2$ we take $r=1$ and so are 
considering the action of $\SU(2)$ on $\C P^1$.)
An important point here is that this Hamiltonian $H_A$ is slow in the 
sense of Section~1.3 so that the path it generates does minimize Hofer length.
This is why the methods sketched below give sharp estimates.

Suppose given a product $\prod_{j=1}^NA_j=\1$ of $N$ natrices.
Entov constructs a 
Hamiltonian bundle  $M\to (E,\Om)\to B$ 
over $B = S^2\setminus\{z_1,\dots,z_N\}$ which is trivial 
topologically but supports a symplectic 
form $\Om$ whose monodromy round the $j$th puncture is (conjugate to) $A_j$.
When $\prod_jA_j=\1$ it is possible to construct $E$ by cutting and 
pasting so that the area of  $(E,\Om)$  is arbitrarily close to zero.
Then one chooses an almost complex structure $\TJ$ on $E$ 
so that the projection to $B$ (with its obvious structure)
is holomorphic and so that $\TJ$ is normalized
 near each puncture to be compatible with the monodromy.  One can do this in 
 such a way that, in the obvious product coordinates  near each end,
 the $\TJ$-holomorphic sections are graphs of solutions of the
 Floer trajectory equation~(\ref{eq:Fl}) for 
 the Hamiltonian functions $H_{j}$ that generate the monodromy $A_j$.
 This has several consequences:
 
 \MS

 \NI
 $\bullet$   at the $j$th end 
 each $\TJ$-holomorphic section $\tu:B\to E$ of finite 
 energy converges to some critical point  $x_{I_j}\in M$
 of the appropriate Hamiltonian $H_{j}$;\MS
 
 \NI
  $\bullet$ the symplectic area 
  $\int_B\tu^*(\Om) $ of the section (which has to be positive)  is 
  now the sum of three terms;
  the area of $(E,\Om)$, the term $\om(A)$ where $A\in H_2(M)$ 
  measures the homology class of $\tu$, and the boundary contribution, which is
  precisely equal to the sum of the critical values $H_{A_j}(x_{I_j})$  at the 
  limiting critical points.  Thus for each section $\tu$ one obtains an 
  inequality
  $$
  \sum_{j=1}^{N} -H_{A_j}(x_{I_j}) \le \om(A).
  $$
  Entov calls these {\bf action inequalities}.  They are generated in 
  the same way as
 the area inequality of Proposition~\ref{prop:main}.  Although they 
 seem rather different, because they come from a section of a trivial 
 bundle over a noncompact space rather than a section of a 
 nontrivial bundle 
 over $S^2$, the next paragraph explains that there is little 
essential  difference in the set up. 
  
  \MS
  
  We need  to understand the moduli spaces of sections of $(E,\Tilde J)\to 
  B$.  When are they nonempty?  Here we are looking at 
  sections over a punctured sphere with Floer boundary conditions.  
  One can \lq\lq cap off'' these ends (by the same gluing arguments 
  that establish the PSS isomorphism between quantum and Floer 
  homology), constructing from such a section $\tu$ a 
  section $\tv$ of 
  the trivial bundle $S^2\times M\to S^2$ such that
  $\tv(z_j)\in z_j\times C_{I_j}$, where $C_{I_j}$ is the 
  Schubert cycle corresponding to the critical point $x_{I_j}$, i.e. 
  the unstable manifold of $x_{I_j}$ under the downward gradient flow 
  of $H_j$.   Conversely, given a $J$-holomorphic sphere $v:S^2 \to M$ 
  such that $v(z_j)\in C_{I_j}$ for all $j$, one can construct from its graph 
  $\tv$ a $\TJ$-holomorphic  
  section $\tu$ of $(E,\TJ)\to B$ with the corresponding limiting behaviour.
  Therefore we need to understand when these spheres $v$ exist.
  Here the points $z_j\in S^2$ are fixed.  Therefore, the number of 
  such maps $v$ is measured by the Gromov--Witten invariant
  \begin{equation}\label{eq:GW}
  \GW^{M,\{1,\dots,N\}}_{A,N}([C_{I_1}],\dots,[C_{I_N}]).
  \end{equation}
  The superscript $\{1,\dots,N\}$ on $\GW$ indicates that the marked points 
  are fixed, so that this invariant measures the number of 
  intersections of  $[C_{I_1}]\times \dots\times [C_{I_N}]$ with the 
  evaluation map 
  $$
  \ev:\Mm(M,A;J) \to M^N,\quad u\mapsto u(z_1,\dots,z_N).
  $$
 These invariants (or correlators)
  are nonzero precisely when the $A$-component of the
 quantum product $[C_{I_1}]*\dots*[C_{I_{N-1}}]$ has nontrivial 
 intersection with the Schubert class  $[C_{I_N}]$.    These are 
  the nontrivial correlators  in quantum Schubert calculus
  to which  Agnihotri--Belkale refer.  Whenever one such correlator is 
  nonzero, the corresponding moduli space of sections  cannot be empty, and hence 
  one gets an inequality which turns out to be precisely one of the ABW
  inequalities.
 \MS
 
 This is just one application of Entov's work.  He generalizes 
 Proposition~\ref{prop:polt} and Corollary~\ref{cor:main}
 to give an interpretation of the 
 action inequalities in terms of an appropriate Hofer distance between the 
 conjugacy classes in $\Ham(M)$ containing the elements $A_i$.

\section{Existence of length minimizing paths in 
$\Ham(M)$.}\label{sec:lect6}

Recall that a Hamiltonian $H_t, t\in [0,1],$ is said to have a {\bf fixed minimum}
at the point $p$ if each function $H_t$ takes its minimum value at $p$.
In this lecture we sketch the proof of the following result from~\cite{McvH}.
Oh gives a quite different proof of this in~\cite{Oh2}.

\begin{thm}\label{thm:lmin}  If $H_t, t\in [0,1],$ has both a fixed 
maximum and a fixed minimum and if it is  sufficiently small in the 
$C^2$-norm, then the path $\phi_t^H$ that it generates in $\Ham(M)$ 
minimizes the Hofer norm, i.e.
$$
\rho(\phi_1^H) = \Ll(\phi_t^H).
$$
\end{thm}

\NI
{\bf Remarks}

\NI(i)
The existence of the fixed extrema is necessary: any path without 
a fixed minimum for example can be altered (keeping the endpoints fixed)
so as to preserve $\Ll^+$ but decrease 
$\Ll^-$.  Therefore for $C^2$ small paths the above proposition gives 
a necessary and sufficient condition for them to realise the Hofer 
norm.  

\NI(ii) It is possible to extract from the proof
 a precise description of how small $H_t$ 
must be for the above result to hold.  
The bound depends only on $(M,\om)$.\MS

\subsection{Idea of the proof}

  If $\rho(\phi_1^H) < \Ll(\phi_t^H)$ then 
there is another shorter Hamiltonian, say $K_t$, with the same time-$1$ map.
Therefore either $\Ll^+(K_t) < \Ll^+(H_t)$ or
$\Ll^-(K_t) < \Ll^-(H_t)$; say the former.  We then consider the space
$(R_{K,H}(2\eps), \Om_0)$ formed by gluing the thickened region 
$(R_H^-(\eps), \Om_0)$ 
under the graph of $H$ to the region $(R_{K}^+(\eps),\Om_0)$ above the graph 
of $K$ along the monodromy of the hypersurfaces $\Ga_K$ and $\Ga_{H}$.
(For definitions, see the proof of Proposition~\ref{prop:polt}.)
This gives rise to a space $(R_{K,H}(2\eps), \Om_0)$ with trivial monodromy round 
its boundary and that fibers over a disc.  Identifying the boundary of this 
disc to a point, one therefore gets a symplectic fibration 
$$
(P_{K,H}(2\eps), \Om_0) \to S^2.
$$
By construction,
$$
\area(P_{K,H}(2\eps), \Om_0)
\;\;=\;\; \Ll^+(K_t) + \Ll^-(H_t) + 2\eps \;\;<\;\; \Ll(H_t), 
$$
provided that $\eps$ is sufficiently small.

Next we use the following fact from~\cite{LM} which is proved by a simple 
geometric construction.
Recall that the capacity of a ball of radius $r$ is $\pi r^{2}$.

\begin{lemma}\label{le:1}
If $H_{t}$ is sufficiently small in the $C^{2}$-norm and has a fixed 
maximum (resp. minimum), then for all 
$\eps > 0$ it is 
possible to embed a ball of capacity $\Ll(H_{t})$ in 
$R_{H}^{-}(\eps)$ (resp. $R_{H}^{+}(\eps)$).
\end{lemma}

Therefore the manifold $(P_{K,H}(2\eps), \Om_0)$ contains an embedded ball
with capacity larger that its area.  If the fibration $(P_{K,H}(2\eps), 
\Om_0)\to S^2$ were symplectically trivial, this would contradict the 
nonsqueezing theorem proven in~\cite{LM} for
so-called \lq\lq quasicylinders".  As it is, we have no control on the topology of this 
fibration:  it is built from the loop $\la = (\phi_{t}^K) 
*(\phi_t^H)^{-1}$
which does not have to contract in $\Ham(M)$. 
Therefore, what we need to do is prove a version of the nonsqueezing 
theorem that holds in this context.  

The rest of the lecture will discuss this question.  Full details of the 
argument outlined above can be found in~\cite{McvH,LM}.

\subsection{Nonsqueezing for fibrations of small area}

\begin{defn} We say that {\bf the nonsqueezing theorem holds for 
the fibration} $(P,\Om)\to S^2$ if area$(P,\Om)$ constrains 
the radius of any embedded  symplectic ball $B^{2n+2}(r)$ 
 in $(P,\Om)$ by the inequality
$$
\pi r^{2} \;\le\; {\rm area\,} (P, \Om).
$$
\end{defn}

Here is a question that is still open for arbitrary manifolds $(M, \om)$.  

\begin{question}  Is there an $\eps = \eps(M,\om) > 0$ such that
the nonsqueezing theorem holds for all fibrations $(P_\la,\Om)\to S^2$
whose generating loop $\la$ has length $\trho^+(\la)\le \eps$?  Would this 
be true if we bounded the length of both sides of $\la$, i.e.
we assumed that both $\trho^+(\la)$ and $\trho^-(\la) = \trho^+(-\la)$
are $\le \eps$?
\end{question}

An affirmative answer (to either question) would be enough to finish the 
proof of Theorem~\ref{thm:lmin}.  For, by choosing $H_t$ so small that 
$\Ll(H_t) <\eps/2$ we could ensure that both $(P_{K,H}(2\de), \Om_0)$
and $(P_{H,K}(2\de), \Om_0)$ had area $< \eps$. But they both contain 
balls of capacity $= \Ll(H_t)$ by Lemma~\ref{le:1}, and one of them has to 
have area $< \Ll(H_t)$. 

We  show  in~\cite{McvH}
that if $(M, \om)$ is a spherically integral symplectic 
manifold (i.e. $[\om] \in H^2(M, \Z)$), the 
nonsqueezing theorem holds for all loops $\la$ in $\Ham(M, \om)$
with $\trho^+(\la) + \trho^-(\la) < 1/2.$  Thus in 
this case we may take $\eps = 1/2$.

The best result for general manifolds involves the idea of {\bf weighted 
nonsqueezing}.  In other words the nonsqueezing inequality must be modified 
by a weight $\ka_0$.  Now the size of $\eps$ is governed by the 
constant $\hbar$ of equation~(\ref{eq:Hbar}).

\begin{prop}\label{prop:ns3}  Suppose that $\la$ is a loop in 
$\pi_1(Ham(M), \om))$ such that $\trho^+(\pm \la) < \hbar(M)/2$. Then there is 
 $\ka_0\in \R$ with $|\ka_0| \le \max(\trho^+(\la), \trho^+(-\la))$
such that the radii of all symplectically 
embedded balls in $(P_{\pm\la}, \Om)$ are 
constrained by the inequalities
$$
\pi r^2 \le \area(P_{\la}, \Om) + \ka_0,\quad 
\pi r^2 \le \area(P_{-\la}, \Om) - \ka_0.
$$
\end{prop}

\NI
{\bf Proof of Theorem~\ref{thm:lmin} assuming 
Proposition~\ref{prop:ns3}.}\MS

Suppose as above that $\phi_1^K = \phi_1^H$ has Hofer norm $<\hbar/4$
and that
 $$
\Ll^{+}(K_{t}) + \Ll^{-}(K_{t}) \;= \;\Ll(\ga) - \de \;<\;  
\Ll^{+}(H_{t}) + \Ll^{-}(H_{t}) < \hbar/4.
$$ 
As before, we may assume that:
$$
\Ll^{+}(K_{t}) \;\; = \;\;  \Ll^{+}(H_{t}) - \de' 
\;\; < \;\; \Ll^{+}(H_{t}),\qquad
\Ll^{-}(K_{t}) \;\; = \;\; \Ll^{-}(H_{t}) - \de + \de',
$$
for some $\de' > 0$.   Let   $\la =\phi_t^K\circ (\phi_t^H)^{-1} $ as before
so that $P_{K,H} = P_\la$, $P_{H,K} = P_{-\la}$.
Then for small $\eps $
$$
\begin{array}{lcl}
\area(P_{K,H}(\eps), \Om_{0}) & = &  \Ll(H_t)  - \de' + \eps \;\;<\;\; \Ll(H_t)
\;\; < \;\;
\hbar/4,\\ 
\area(P_{H,K}(\eps), \Om_0) & = & \Ll(H_t) - \de + \de' + \eps \;\;\le\;\;  
2\Ll(H_t) < \hbar/2.
\end{array}
$$
By Proposition~\ref{prop:ns3} there is $\ka_0$ with $|\ka_0| \le \hbar/2$
 such that embedded balls satisfy
$$
\pi r^2 \le \area(P_{K,H}(\eps), \Om) + \ka_0,\quad 
\pi r^2 \le \area(P_{H,K}(\eps), \Om) - \ka_0.
$$
But, by construction, both $(P_{K,H}(\eps), \Om)$ and $(P_{H,K}(\eps), \Om)$ 
contain embedded balls of 
capacity $\pi r^2 = \Ll(\ga) > \area(P_\la(\ka_0), \Om)$.  Hence 
 $\ka_0 > 0$.  Further,
$$
\begin{array}{lclcl}
\Ll(H_t) &\le & \area(P_{K,H}(\eps), \Om) + \ka_0 & = &\Ll(H_t) -\de' + 
\eps + \ka_0\\
\Ll(H_t) & \le & \area(P_{H,K}(\eps), \Om) - \ka_0 & = &
\Ll(H_t) - \de + \de' + \eps - \ka_0.
\end{array}
$$
Adding, we find $0 \le -\de + 2\eps$.  Since $\de$ is positive and 
 $\eps$ can be arbitrarily small, this 
is impossible.  Hence result.\QED

It therefore remains to prove Proposition~\ref{prop:ns3}. Recall the 
definition of $\Ss(\la)$ from~(\ref{eq:qm2})
It will be convenient to separate out the terms in $\Ss(\la)$
with $\al_{d,\ka} \in H_{2n}(M)$. Since $\Ss(\la)$ has degree $2n$,
such terms must have $d=0$.  Thus we write
$$
\Ss(\la) = \sum_{0,\ka}\ r_{\ka}\1\otimes t^\ka + x',
$$
where $\1 = [M]$ is the unit element, $r_\ka\in \Q$ and
$$
x' \in \QH^+ : = \sum_{i>0}\; H_{2n-i}(M)\otimes \La.
$$

\begin{defn}\label{def:sns} We  say that the fibration $(P, 
\Om)\to S^{2}$  with fiber $M$ has  a
{\bf good section of weight} $\ka_0$ if there is 
a class $\TB \in H_{2}(P)$ such that
\begin{itemize}
\item[(i)] $ \GW^{P}_{\TB,3}([M], [M], pt) \ne 0$;
\item[(ii)]  $u(\TB) = - \ka_0$ where $u$ is the coupling class.
\end{itemize}
Note that $\ka_0$ could be positive or negative.
\end{defn}

In particular,  in the language of Lemma~\ref{le:hbar}, we can take
$$
\ka_0 = v(\sum r_{\ka}\1\otimes t^\ka).
$$

\begin{lemma}\label{le:ns1} 
Suppose that $(P_{\la},\Om)$ has a good section of weight $\ka_0$.
  Then the radius $r$ of an 
embedded ball in $(P,\Om)$ is constrained by the inequality:
$$
\pi r^2 \le \area\,(P, \Om) - \ka_0.
$$
\end{lemma}
\proof{}  The hypotheses imply that for some section class $\TB$ 
with $u_\la(\TB) = -\ka_0$ we have
 $$
 \GW^P_{\TB,3}([M], [M], pt) = r_\ka\ne 0.
 $$
Since this invariant counts perturbed $J$-holomorphic stable maps in class 
$\TB$ through an arbitrary point, it follows that there is such 
a curve through every point in $P$.  Since the perturbation can be taken 
arbitrarily small, it follows from Gromov's compactness theorem that there 
has to be some $J$-holomorphic stable map in this class through every point 
in $P$.
 Hence the usual arguments (cf.~\cite{GRO} or~\cite{LMe})
 imply that the radius  $r$ of any embedded ball satisfies the 
 inequality:
 $$
 \pi r^{2} \;\;\le\;\; [\Om](\TB) \;\; \le 
\;\; \area(P_{\la}, \Om) - \ka_0,
 $$
where the last inequality follows as in Proposition~\ref{prop:main}.
The result follows.\QED

\NI
{\bf Proof of Proposition~\ref{prop:ns3}}\MS

By hypothesis there are fibrations $(P_{\pm\la}, \Om)$ each with area $< 
\hbar/2$.  Choose $\de > 0$ so that
\begin{equation}\label{eq:we}
\trho^{+}(\la) + \trho^{-}(\la) < \hbar - 2\de.
\end{equation}
By  Proposition~\ref{prop:polt}, there is a $\om$-compatible symplectic form
$\Om_{\la}$ on $P_{\la}$ with area $< \trho^{+}(\la) + \de$, and  a 
similar form $\Om_{-\la}$ on $P_{-\la}$ with area $< \trho^{+}(-\la) + \de$.
Write
$$
\Ss(\la) =  \sum_{\ka} \;\1\otimes\,r_{\ka} t^{\ka} + x,
\quad
\Ss(-\la) =  \sum_{\ka}\; \1\otimes \,r_{\ka}' t^{\ka} + x'
$$
where $r_{\ka},r_{\ka}'\in \Q$ and $x,x'\in \QH^{+}$.
 Proposition~\ref{prop:main} implies that
%  if $r_\ka\ne 0$ then
% 
%  for all  $\ka$  that occur in $\Ss(\la)$ 
% (resp.  $\Ss(-\la)$) with
% nonzero coefficient
$$
r_\ka\ne 0\;\;\Longrightarrow \ka\; < \;\;\trho^{+}(\la) + \de,\qquad 
r_\ka'\ne 0\;\;\Longrightarrow \ka\; < \;\;\trho^{+}(-\la) + \de
$$

Next apply the valuation $v$ in~(\ref{eq:val}) to the identity
$$
\Ss(\la)*\Ss(-\la) = \Ss(0) = \1.
$$
We claim that at least one of $\Ss(\la),  \Ss(-\la)$
has a nonzero term $\1\otimes r_\ka t^{\ka}$ with $\ka\ge 0$.
%  with $q_{B}\ne 0, \om(B) \ge 0$.
For otherwise since $\Ss(\la)*\Ss(-\la)=\1$,
the product $x*x'$ must contain the term  $\1\otimes 
t^{0}$ with a nonzero 
coefficient.  Because this term appears in $x*x' - x\cap x'$ we find
from Lemma~\ref{le:hbar} that
$$
0 \;=\; v(\1\otimes t^{0}) \;\;\le\;\; v(\Ss(\la)) + v(\Ss(-\la)) - 
\hbar(M)\;\; \le\;\; 
\trho^{+}(\la) + \trho^{+}(-\la) + 2\de  - \hbar \;\;<\;\; 0,
$$
contradicting~(\ref{eq:we}).
 Therefore for $\la'$ equal to at least one of $\la$ or $-\la$, 
$\Ss(\la')$
has a term $r_{\ka}\1\otimes  t^\ka$
 with $r_{\ka}\ne 0$ and $0 \le \ka < \area(P_{\la'}, \Om)$.  

For $\la' = \pm\la$ set 
$$
% \ka(\la') = \max \{\ka: r_\ka\ne 0 \;\mbox{ in }\;\Ss(\la')\}.
$$
 The equation $\Ss(\la)*\Ss(-\la) = \1$ implies that
$\ka(\la) = -\ka(-\la)$.  Moreover,  
 by Lemma~\ref{le:ns1}, the radius $r$ of any embedded ball in 
 $(P_{\la'}, \Om)$ 
satisfies
$$
\pi r^2 \le \area(P_{\la'}, \Om) - \ka(\la').
$$
Hence we may take $\ka_0 = -\ka(\la)$.\QED

This completes the proof of Theorem~\ref{thm:lmin}.

\end{document}